\documentstyle[11pt]{article}

       \font\tenmsb=msbm10
       \font\sevenmsb=msbm7
       \font\fivemsb=msbm5
       \catcode`\@=11
       \ifx\amstexloaded@\relax\catcode`\@=\active
       \endinput\else\let\amstexloaded@\relax\fi
       \def\spaces@{\space\space\space\space\space}
       \def\spaces@@{\spaces@\spaces@\spaces@\spaces@\spaces@}
       \def\space@.{\futurelet\space@\relax}
       \space@. %
       \def\Err@#1{\errhelp\defaulthelp@\errmessage{AmS-TeX error: #1}}
       \def\relaxnext@{\let\next\relax}
       \def\accentfam@{7}
       
       \def\noaccents@{\def\accentfam@{0}}
       \def\Cal{\relaxnext@\ifmmode\let\next\Cal@\else
       \def\next{\Err@{Use \string\Cal\space only in math mode}}\fi\next}
       \def\Cal@#1{{\Cal@@{#1}}}
       \def\Cal@@#1{\noaccents@\fam\tw@#1}

       \def\Bbb{\relaxnext@\ifmmode\let\next\Bbb@\else
       \def\next{\Err@{Use \string\Bbb\space only in math mode}}\fi\next}
       
       \def\Bbb@#1{{\Bbb@@{#1}}}
       \def\Bbb@@#1{\noaccents@\fam\msbfam#1}
       \newfam\msbfam
       \textfont\msbfam=\tenmsb
       \scriptfont\msbfam=\sevenmsb
       \scriptscriptfont\msbfam=\fivemsb

\def\N{{\Bbb N}}
\def\Z{{\Bbb Z}}
\def\R{{\Bbb R}}
\def\T{{\Bbb T}}

\newtheorem{Theorem}{Theorem}
\newtheorem{Lemma}{Lemma}[section]

\@addtoreset{equation}{section}

\newcommand{\beq}{\begin{equation} }
\newcommand{\eeq}{\end{equation} }

\makeatletter

\newcommand{\Rmnum}[1]{\expandafter\@slowromancap\romannumeral #1@}
\makeatother

\begin{document}

\setlength{\columnsep}{5pt}

\title{An abstract infinite dimensional KAM theorem with application to nonlinear higher dimensional Schr\"odinger equation systems \thanks{This work is partially supported
by NSFC grant 11271180.}}
\author{\\Shidi Zhou\\
 {\footnotesize Department of Mathematics,Nanjing University, Nanjing 210093, P.R.China}\\
{\footnotesize Email: mathsdzhou@126.com}}


\date{}
\maketitle

\begin{abstract}
 In this paper we consider nonlinear Schr\"odinger systems with periodic boundary condition in high dimension. We establish an abstract infinite dimensional KAM theorem and apply it to the nonlinear Schr\"odinger equation systems with real Fourier Multiplier. By establishing a block-diagonal normal form, We prove the existence of a class of Whitney smooth small amplitude quasi-periodic solutions corresponding to finite dimensional invariant tori of an associated infinite dimensional dynamical system.
\end{abstract}



\section{\textbf{Introduction}}
 \noindent

 In this paper, we consider a class of higher dimensional nonlinear Schr\"odinger equation systems with real Fourier Multiplier under periodic boundary condition:
\begin{eqnarray}
\left\{
\begin{array}{l}
-i\partial_t u = -\Delta_x u + M_{\xi} u + \partial_{\bar u} G(|u|^2, |v|^2)  \\
-i\partial_t v = -\Delta_x v + M_{\sigma} v + \partial_{\bar v} G(|u|^2, |v|^2)  \\
\end{array}
\right.
\end{eqnarray}
where
$
u=u(t,x),v=v(t,x), t\in \R, x\in \T^d, d\geq 2
$.
$G=G(a,b)$ is a real analytic function defined in a neighbourhood of the origin in $\R^2$ with
$G(0,0)=\partial_a G(0,0) = \partial_b G(0,0)=0$,
 which means that the Taylor series of $G$ with respect to $(a,b)$ should start from the second order term (a classical example of the nonlinearity is $|u|^2 |v|^2$ ). $M_{\xi},M_{\sigma}$ are two real Fourier Multipliers, which supply artificial parameters (defined in section 2). By a KAM algorithm, we prove that after removing a small measure of parameters, the equation system (1.1) admits a class of small-amplitude quasi-periodic solutions under sufficiently small nonlinear perturbations.

KAM theory has been a fundamental tool for years in studying Hamiltonian PDEs by constructing a class of invariant torus. Starting from the pioneering work $\cite {22}$ by Kuksin in 1987, Newtonian scheme was developed by $\cite{7,22,34}$ and has showed its power in doing Hamiltonian PDEs in one-dimensional.  The main idea is to construct a local normal form of the solution of the equation, which is fundmental in the study of the dynamical properties of the initial equation, and then carry out an infinite iteration process. Thus a class of invariant torus and corresponding quasi-periodic solutions are got. For related works, see $\cite{6,13,14,20,21,23,24,25,26,27,28,35}$.

Hamiltonian PDEs in high dimension have also attracted great interests. The KAM method working on one-dimensional PDEs are not effective enough here, due to the multiplicity of eigenvalues of the linear operator, which causes terrible resonance between two eigenvalues $\lambda_n=|n|^2+o(1)$ and $\lambda_m=|m|^2+o(1)$ if $|n|=|m|$ ($|\cdot|$ means $l^2$ norm here). This is a big obstacle in realizing the second Melnikov condition at each KAM iterative step. The first breakthrough comes from Bourgain's work $\cite{4}$ in 1998, in which a class of two-frequency quasi-periodic solution of two-dimensional nonlinear Schr\"odinger equations was got. In this paper, Bourgain introduced the famous multi-scale analysis method, which avoided the cumbersome second Melnikov condition. Later in $\cite{5}$, Bourgain improved his method and got the small-amplitude quasi-periodic solutions of high dimensional Schr\"odinger equations and wave equations. Following his idea and method, abundant works have been done, see $\cite{1,2,3,33}$.

Although multi-scale analysis has great advantages, its drawbacks couldn't be ignored. For example, we can't see the local Birkhoff normal form of the equation, and the linear stability of the solution is also unavailable. Thus a KAM approach is also expected in dealing with high-dimensional Hamiltonian PDEs. The first work comes from Geng and You $\cite{15}$ in 2006, in which the quasi-periodic solutions of beam equation and nonlocal smooth Schr\"odinger equation were got. In \cite{15} the nonlinearity of the equation should be independent of the spatial variable $x$, which implies that the Hamiltonian satisfies the important property ``zero-momentum'' (condition $\bf (A4)$ in $\cite{15}$). The multiple eigenvalues of linear operator are avoided by making use of "zero-momentum" condition and measure estimate is conducted by the help of regularity of equation. A much more difficult question: nonlinear Schr\"odinger equation defined on $\T^d,d\geq 2$, with convolutional type potential and nonlinearity dependent on spatial variable $x$, was solve by Eliasson and Kuksin $\cite{9}$ in 2010. In this milestone-style paper, to work with the multiple eigenvalues, they studied the elaborate distribution of integers points on a sphere and according to this they got the normal form of block-diagonal, with each block becoming larger and larger along the KAM iteration. Besides, they developed a very important property ``Lipschitz-Domain'' to do the measure estimate, due to the absence of regularity of the equation. Following their idea and method, Geng , Xu and You $\cite{12}$ got the quasi-periodic solutions of completely resonant Schr\"odinger equation on $\T^2$, by constructing some appropriate tangential sites on $\Z^2$. Later in $\cite{29,30}$, C.Procesi and M.Procesi got the same result of $\cite{12}$ in arbitrary dimension. In $\cite{29,30}$ their had a very ingenious choice of tangential sites through the method of graph theory. For other works on PDEs in high dimension, see $\cite{8,10,11,16,17,18,31,32}$.

Although there has been rich work about Hamiltonian equations, until now little is known about Hamiltonian equation sets, i.e. two coupled equations. In $\cite{19}$ Grebert, Paturel and Thomann contructed the beating solutions of Schr\"odinger equation system in one-dimension and got the growth of Sobolev norms of the solutions. However, quasi-periodic solutions corresponding to finite dimensional invariant torus is still unknown. Our present paper is working on the nonlinear Schr\"odinger system with real Fourier Multiplier. A more interesting question is about the completely resonant equation system, i.e. no artificial parameters are imposed :
\begin{eqnarray}
\left\{
\begin{array}{l}
-i\partial_t u = -\Delta_x u + |v|^2 u  \\
-i\partial_t v = -\Delta_x v + |u|^2 v
\end{array}
\right.
\end{eqnarray}
This equation system will be dealt with in our forthcoming paper.

$\\$
\indent Now let's state our main theorem:

\begin{Theorem}
Let
$
S=\{i_1,\cdots,i_b\} \subseteq \Z^d, \tilde{S}=\{t_1,\cdots,t_{\tilde{b}}\} \subseteq \Z^d, 0\in S \bigcap \tilde{S}; d\geq 2; b,\tilde{b} \geq 2;
$.
There exists a Cantor set ${\mathcal C} \subseteq \R^{b+\tilde{b}}$ of positive measure, s.t.
$\forall (\xi,\sigma) \in {\mathcal C}$, the nonlinear Schr\"odinger equation system (1.1) admits a class of small amplitude quasi-periodic solutions of the form:
\begin{eqnarray}
\left\{
\begin{array}{l}
u(t,x)=\sum\limits_{j=1} ^b c_j e^{i\omega_j t} \phi_{i_j} + O(|\xi|^{\frac{3}{2}}+|\sigma|^{\frac{3}{2}})
\quad \omega_j = |i_j|^2 + O(|\xi|+|\sigma|) \\
v(t,x)=\sum\limits_{j=1} ^{\tilde{b}} \tilde{c}_j e^{i\tilde{\omega}_j t} \phi_{t_j} + O(|\xi|^{\frac{3}{2}}+|\sigma|^{\frac{3}{2}})
\quad  \tilde{\omega}_j = |t_j|^2 + O(|\xi|+|\sigma|)
\end{array}
\right.
\end{eqnarray}
\end{Theorem}

The rest of the paper is organized as follows: In section 2 we introduce some notations and state the abstract KAM theorem; In section 3 we deal with the normal form; In section 4, we conduct one step of KAM iteration; In section 5 we state the iterative lemma; In section 6 we do the measure estimate.

\section{\textbf{Preliminaries and statement of the abstract KAM theorem}}
\noindent

 In this section we introduce some notations and state the abstract KAM theorem. The KAM theorem can be applied to (1.1) to prove Theorem 1.

Given two set $S,\tilde S \subseteq \Z^d, d \geq 2$, $S=\{i_1,i_2,\cdots,i_b\},
\tilde S=\{t_1,t_2,\cdots,t_{\tilde b}\}, b,\tilde b \geq 2$ (for convenience we assume $0\in S \bigcap \tilde S$).
 Let $\Z^d _1$ be the complementary set of $S$ in $\Z^2$, and $\Z^d _2$ be the complementary set of $\tilde S$ in $\Z^2$. Denote $u=(u_n)_{n\in \Z^d _1}$ with its conjugate $\bar u=(\bar u_n)_{n\in \Z^d _1}$, and similarly $v=(v_n)_{n\in \Z^d _2}$ with its conjugate $\bar v=(\bar v_n)_{n\in \Z^d _2}$.
We introduce the weighted norm as follows:
\begin{eqnarray}
\| u \|_{\rho} = \sum_{n\in \Z^d _1} |u_n| e^{\rho |n|},
\| v \|_{\rho} = \sum_{n\in \Z^d _2} |v_n| e^{\rho |n|} \qquad \rho >0
\end{eqnarray}
\noindent here $|n|=\sqrt{|n_1|^2 +\cdots + |n_d|^2}$, $n=(n_1,\cdots,n_d)\in \Z^d _1$ for $u$ and in $\Z^d _2$ for $v$ . Denote a neighborhood of $\T^{b+\tilde b}\times (\{I=0\}\times \{J=0\}) \times
(\{u=0\}\times\{\bar u=0\})\times (\{v=0\}\times\{\bar v=0\})$ by
\begin{eqnarray}
D(r,s)=\{
(\theta,\varphi, I,J, u, \bar u,v,\bar v): |{\rm Im} \theta|, |{\rm Im} \varphi| < r;
I, J< s^2;  \| u \|_{\rho},\| \bar u \|_{\rho}, \| v \|_{\rho},\| \bar v \|_{\rho} < s
\}
\end{eqnarray}
\noindent where $|\cdotp|$ means the sup-norm of complex vectors.

Let $\alpha = \{\alpha_n\}_{n\in \Z^d _1}, \beta = \{\beta_n\}_{n\in
\Z^d _1},\tilde{\alpha} = \{\tilde{\alpha}_n\}_{n\in \Z^d _2}, \tilde{\beta} = \{\tilde{\beta}_n\}_{n\in
\Z^d _2}$, $\alpha_n, \beta_n, \tilde{\alpha}_n, \tilde{\beta_n} \in \N$ with only finitely many
non-vanishing  components. Denote $u^{\alpha}\bar
u^{\beta}=\prod\limits_{n\in \Z^d _1}u^{\alpha_n}_n \bar u^{\beta_n}_n,
v^{\tilde{\alpha}}\bar v^{\tilde{\beta}}=\prod\limits_{n\in \Z^d _2}v^{\tilde{\alpha}_n}_n \bar v^{\tilde{\beta}_n}_n$
 and let
\begin{eqnarray}
F(\theta,\varphi, I,J, u,\bar u,v,\bar v)=\sum_{k l \alpha \beta,\tilde{k} \tilde{l} \tilde{\alpha} \tilde{\beta}}
F_{kl\alpha\beta,\tilde{k} \tilde{l} \tilde{\alpha} \tilde{\beta}}(\xi,\sigma)
e^{i(\langle k,\theta \rangle+\langle \tilde{k},\varphi \rangle)}I^l J^{\tilde{ l}}
u^\alpha \bar u^\beta v^{\tilde{\alpha}} \bar v^{\tilde{\beta}}
\end{eqnarray}
\noindent where $(\xi,\sigma) \in {\mathcal O}\subseteq \R^{b+\tilde b}$ is the
parameter set. $k=(k_1,\cdots,k_b)\in \Z^b, \tilde{k}=(\tilde{k}_1,\cdots,\tilde{k}_{\tilde{b}})\in \Z^{\tilde b}$ and
$l=(l_1,\cdots,l_b)\in \N^b, \tilde{l}=(\tilde{l}_1,\cdots,\tilde{l}_{\tilde b})\in \N^{\tilde b}$,
$I^l = I_1 ^{l_1} \cdots I_b ^{l_b}$,
$J^{\tilde l} = J^{\tilde l_1} _1\cdots J^{\tilde{l} _{\tilde b}}_{\tilde b}$.
Denote the weighted norm of $F$ by
\begin{eqnarray}
&\| F \|_{D(r,s), \mathcal{O}} =
\sup\limits_{(\xi,\sigma) \in {\mathcal O}, \| u \|_{\rho},\| \bar u \|_{\rho},
\| v \|_{\rho},\| \bar v \|_{\rho} < s
}
\sum_{kl\alpha\beta,\tilde{k}\tilde{l}\tilde{\alpha}\tilde{\beta}}
|F_{kl\alpha\beta,\tilde{k}\tilde{l}\tilde{\alpha}\tilde{\beta}}|_{\mathcal{O}}
e^{(|k|+|\tilde k|)r} \nonumber \\
& \times s^{2(|l|+|\tilde l|)}
|u^{\alpha}||\bar u^{\beta}| |v^{\tilde{\alpha}}||\bar v^{\tilde{\beta}}| \\
&|F_{kl\alpha\beta,\tilde k \tilde l \tilde \alpha \tilde \beta}|_{\mathcal{O}}=\sup\limits_{(\xi,\sigma) \in \mathcal{O}}\sum_{0\leq d \leq 4}|\partial_{(\xi,\sigma)}^4 F_{kl\alpha\beta,\tilde k \tilde l \tilde \alpha \tilde \beta}|
\end{eqnarray}
\noindent where the derivatives with respect to $(\xi,\sigma)$ are in the
sense of Whitney.

To a function $F$ we define its Hamiltonian vector field by
\begin{eqnarray}
X_F = (
F_I,F_J, -F_{\theta}, -F_{\varphi}, i\{F_{u_n}\}_{n\in \Z^d _1}, -i\{F_{\bar u_n}\}_{n\in \Z^d _1},
i\{F_{v_n}\}_{n\in \Z^d _2}, -i\{F_{\bar v_n}\}_{n\in \Z^d _2}
)
\end{eqnarray}
and the associated weighted norm is
\begin{eqnarray}
\|X_F\|_{D(r,s),{\mathcal O}} &=&
\|F_I\|_{D(r,s),{\mathcal O}} + \|F_J\|_{D(r,s),{\mathcal O}} +
\frac{1}{s^2}\bigg(\|F_{\theta}\|_{D(r,s),{\mathcal O}}+\|F_{\varphi}\|_{D(r,s),{\mathcal O}} \bigg)\nonumber \\
&+& \frac{1}{s}
 \left(
 \sum_{n\in \Z^d _1} \bigg(\|F_{u_n}\|_{D(r,s),{\mathcal O}}+\|F_{\bar u_n}\|_{D(r,s),{\mathcal O}}\bigg)e^{|n|\rho}
   \right) \nonumber \\
 &+& \frac{1}{s}
 \left(
 \sum_{n\in \Z^d _2} \bigg(\|F_{v_n}\|_{D(r,s),{\mathcal O}}+
\|F_{\bar v_n}\|_{D(r,s),{\mathcal O}}\bigg)e^{|n|\rho}
 \right)
\end{eqnarray}
\noindent where $\rho > 0$ is a constant and it will shrink at each
iterative step to make the small divisor condition hold due to the lack of regularity of the equation.

$\\$ \indent The normal form $H_0 = N+{\mathcal B}$ with
\begin{eqnarray}
N&=&\langle \omega(\xi,\sigma),I\rangle+\langle \tilde{\omega}(\xi,\sigma),J\rangle +
  \sum_{n\in \Z^d _1}\Omega_n(\xi,\sigma) u_n \bar u_n +\sum_{n\in \Z^d _2}\tilde{\Omega}(\xi,\sigma)_n
  v_n \bar v_n  \\
{\mathcal B}&=&\sum_{n \in \Z^d _1 \cap \Z^d _2}(a_n(\xi,\sigma) u_n \bar v_n + b_n(\xi,\sigma) \bar u_n v_n)
\end{eqnarray}
where $(\xi,\sigma) \in \mathcal O$ is the parameter. Notice that apart from integrable terms, $u_n$ and $\bar v_n$, $v_n$ and $\bar u_n$ may also be coupled and as a result our normal form is in the form of block-diagonal with each block of degree 2.

For this unperturbed system, it's easy to see that it admits a
special solution
$$
(\theta,\varphi,0,0,0,0,0,0)\rightarrow (\theta+\omega t,\varphi+\tilde{\omega} t,
0,0,0,0,0,0)
 $$
 corresponding to an invariant torus in the phase space. Our
goal is to prove that, after removing some parameters, the perturbed
system $H=H_0 + P$ still admits invariant torus provided that $
\|X_P\|_{D_{\rho} (r,s), {\mathcal O}} $ is sufficiently small. To
achieve this goal, we require that the Hamiltonian $H$ satisfies some
conditions:

$\\$ $(\bf A1)$ Nondegeneracy: The map
$(\xi,\sigma) \rightarrow
(\omega (\xi,\sigma),\tilde{\omega} (\xi,\sigma))$
is a $C^4 _W$ diffeomorphism between $\mathcal O$ and its
image ($C^4 _W$ means $C^4$ in the sense of Whitney).

\noindent $(\bf A2)$ Asymptotics of normal frequencies:
\begin{eqnarray}
&\Omega_n = |n|^2 + \acute{\Omega}_n, \quad n\in \Z^d _1 \\
&\tilde{\Omega}_n = |n|^2 + \acute{\tilde{\Omega}}_n, \quad n\in \Z^d _2
\end{eqnarray}
here $\acute{{\Omega}}_n, \acute{\tilde{{\Omega}}}_n$ are $C^4 _W$ functions of $(\xi,\sigma)$.

\noindent $(\bf A3)$ Melnikov conditions: Let
\begin{eqnarray}
A_n =
\left(
\begin{array}{cc}
\Omega_n  & a_n \\
b_n & \tilde{\Omega}_n
\end{array}
\right)  \qquad n\in \Z^d _1 \cap \Z^d _2
\end{eqnarray}
and
\begin{eqnarray}
&A_n = \Omega_n \qquad n\in \Z^d _1 \setminus \Z^d _2   \\
&A_n = \tilde{\Omega}_n \qquad n\in \Z^d _2 \setminus \Z^d _1
\end{eqnarray}
There exists $\gamma,\tau > 0$, s.t. for any $|k|+|\tilde{k}| \leq K, n\in \Z^d _1, m \in \Z^d _2$, one has
\begin{eqnarray}
| \langle k,\omega \rangle + \langle \tilde{k},\tilde{\omega} \rangle | \geq \frac{\gamma}{K^{\tau}}, \qquad
|k|+|\tilde{k}| \neq 0
\end{eqnarray}
and
\begin{eqnarray}
|\det\bigg(
(\langle k,\omega \rangle + \langle \tilde{k},\tilde{\omega} \rangle) I + A_n \bigg)
|\geq \frac{\gamma}{K^{\tau}}
\end{eqnarray}
and
\begin{eqnarray}
|\det \bigg(
(\langle k,\omega \rangle+\langle \tilde{k},\tilde{\omega} \rangle) I
\pm A_n \otimes I_2 \pm I_2 \otimes {A_m}^{\mathrm T}
\bigg)| \geq \frac{\gamma}{K^{\tau}}
\end{eqnarray}
Here $A^{\mathrm T}$ denotes the transpose of matrix $A$ and $I$ denotes the identity matrix.

\noindent $(\bf A4)$ Regularity:
${\mathcal B} + P$ is real analytic with respect to $\theta,\varphi,I,J,u,\bar u,v,\bar v$ and Whitney smooth with respect to $(\xi,\sigma)$. And we have
\begin{eqnarray}
\|X_{\mathcal B}\|_{D_{\rho}(r,s),{\mathcal O}} < 1, \qquad
\|X_P\|_{D_{\rho}(r,s),{\mathcal O}}<\varepsilon
\end{eqnarray}

\noindent $(\bf A5)$ Zero-momentum condition:
The normal form part $ {\mathcal B} + P$ belongs to a class of functions $\mathcal A$ defined by:
$$
f= \sum\limits_{k
\in \Z^b,\tilde{k} \in \Z^{\tilde{b}}, l \in \N^b,\tilde{l} \in \N^{\tilde{b}}, \alpha, \beta,
\tilde{\alpha}, \tilde{\beta}} f_{kl\alpha \beta,\tilde{k}\tilde{l}\tilde{\alpha} \tilde{\beta}}
e^{i(\left\langle k,\theta \right\rangle)+(\langle \tilde{k},\varphi\rangle)}I^l J^{\tilde{l}}
u^{\alpha} {\bar u}^{\beta} v^{\tilde{\alpha}} {\bar v}^{\tilde{\beta}}, \qquad f\in \mathcal A
$$
implies
$$
f_{kl\alpha \beta,\tilde{k}\tilde{l}\tilde{\alpha} \tilde{\beta}}
\neq 0 \Longrightarrow
\sum_{j=1} ^b k_j i_j + \sum_{j=1} ^{\tilde{b}} \tilde{k}_j t_j +
\sum\limits_{n\in \Z^d _1}(\alpha_n - \beta_n)n + \sum\limits_{n\in \Z^d _2}
(\tilde{\alpha}_n - \tilde{\beta}_n)n = 0
$$

\noindent $(\bf A6)$ T\"oplitz-Lipschitz property: There exists a $K>0$.

$\bf (1)$ For any fixed $n,m \in \Z^d _1, c\in \Z^d \setminus \{0\}$, the limits
\begin{eqnarray}
\lim_{t\rightarrow \infty}\frac{\partial^2 P}{\partial u_{n+tc} \partial u_{m-tc}},
\lim_{t\rightarrow \infty}\frac{\partial^2 P}{\partial u_{n+tc} \partial \bar u_{m+tc}},
\lim_{t\rightarrow \infty}\frac{\partial^2 P}{\partial \bar u_{n+tc} \partial \bar u_{m-tc}}
\end{eqnarray}
exists, and moreover, when $|t|>K$, $P$ satisfies:
\begin{eqnarray}
&\|
\frac{\partial^2 P}{\partial u_{n+tc} \partial u_{m-tc}} -
\lim\limits_{t\rightarrow \infty}\frac{\partial^2 P}{\partial u_{n+tc} \partial u_{m-tc}}
\|_{D_{\rho}(r,s),{\mathcal O}} \leq
\frac{\varepsilon}{|t|}e^{-|n+m|\rho}  \\
&\|
\frac{\partial^2 P}{\partial u_{n+tc} \partial \bar u_{m+tc}} -
\lim\limits_{t\rightarrow \infty}\frac{\partial^2 P}{\partial u_{n+tc} \partial \bar u_{m+tc}}
\|_{D_{\rho}(r,s),{\mathcal O}} \leq
\frac{\varepsilon}{|t|}e^{-|n-m|\rho}  \\
&\|
\frac{\partial^2 P}{\partial \bar u_{n+tc} \partial \bar u_{m-tc}} -
\lim\limits_{t\rightarrow \infty}\frac{\partial^2 P}{\partial \bar u_{n+tc} \partial \bar u_{m-tc}}
\|_{D_{\rho}(r,s),{\mathcal O}} \leq
\frac{\varepsilon}{|t|}e^{-|n+m|\rho}
\end{eqnarray}

$\bf (2)$ For any fixed $n\in \Z^d _1, m \in \Z^d _2, c\in \Z^d \setminus \{0\}$, the limits
\begin{eqnarray}
\lim_{t\rightarrow \infty}\frac{\partial^2 P}{\partial u_{n+tc} \partial v_{m-tc}},
\lim_{t\rightarrow \infty}\frac{\partial^2 P}{\partial u_{n+tc} \partial \bar v_{m+tc}},
\lim_{t\rightarrow \infty}\frac{\partial^2 P}{\partial \bar u_{n+tc} \partial \bar v_{m-tc}}
\end{eqnarray}
exists, and moreover, when $|t|>K$, $P$ satisfies:
\begin{eqnarray}
&\|
\frac{\partial^2 P}{\partial u_{n+tc} \partial v_{m-tc}} -
\lim\limits_{t\rightarrow \infty}\frac{\partial^2 P}{\partial u_{n+tc} \partial v_{m-tc}}
\|_{D_{\rho}(r,s),{\mathcal O}} \leq
\frac{\varepsilon}{|t|}e^{-|n+m|\rho}  \\
&\|
\frac{\partial^2 P}{\partial u_{n+tc} \partial \bar v_{m+tc}} -
\lim\limits_{t\rightarrow \infty}\frac{\partial^2 P}{\partial u_{n+tc} \partial \bar v_{m+tc}}
\|_{D_{\rho}(r,s),{\mathcal O}} \leq
\frac{\varepsilon}{|t|}e^{-|n-m|\rho}  \\
&\|
\frac{\partial^2 P}{\partial \bar u_{n+tc} \partial v_{m+tc}} -
\lim\limits_{t\rightarrow \infty}\frac{\partial^2 P}{\partial \bar u_{n+tc} \partial v_{m+tc}}
\|_{D_{\rho}(r,s),{\mathcal O}} \leq
\frac{\varepsilon}{|t|}e^{-|n-m|\rho}  \\
&\|
\frac{\partial^2 P}{\partial \bar u_{n+tc} \partial \bar v_{m-tc}} -
\lim\limits_{t\rightarrow \infty}\frac{\partial^2 P}{\partial \bar u_{n+tc} \partial \bar v_{m-tc}}
\|_{D_{\rho}(r,s),{\mathcal O}} \leq
\frac{\varepsilon}{|t|}e^{-|n+m|\rho}
\end{eqnarray}

$\bf (3)$ For any fixed $n,m \in \Z^d _2, c\in \Z^d \setminus \{0\}$, the limits
\begin{eqnarray}
\lim_{t\rightarrow \infty}\frac{\partial^2 P}{\partial v_{n+tc} \partial v_{m-tc}},
\lim_{t\rightarrow \infty}\frac{\partial^2 P}{\partial v_{n+tc} \partial \bar v_{m+tc}},
\lim_{t\rightarrow \infty}\frac{\partial^2 P}{\partial \bar v_{n+tc} \partial \bar v_{m-tc}}
\end{eqnarray}
exists, and moreover, when $|t|>K$, $P$ satisfies:
\begin{eqnarray}
&\|
\frac{\partial^2 P}{\partial v_{n+tc} \partial v_{m-tc}} -
\lim\limits_{t\rightarrow \infty}\frac{\partial^2 P}{\partial v_{n+tc} \partial v_{m-tc}}
\|_{D_{\rho}(r,s),{\mathcal O}} \leq
\frac{\varepsilon}{|t|}e^{-|n+m|\rho}  \\
&\|
\frac{\partial^2 P}{\partial v_{n+tc} \partial \bar v_{m+tc}} -
\lim\limits_{t\rightarrow \infty}\frac{\partial^2 P}{\partial v_{n+tc} \partial \bar v_{m+tc}}
\|_{D_{\rho}(r,s),{\mathcal O}} \leq
\frac{\varepsilon}{|t|}e^{-|n-m|\rho}  \\
&\|
\frac{\partial^2 P}{\partial \bar v_{n+tc} \partial \bar v_{m-tc}} -
\lim\limits_{t\rightarrow \infty}\frac{\partial^2 P}{\partial \bar v_{n+tc} \partial \bar v_{m-tc}}
\|_{D_{\rho}(r,s),{\mathcal O}} \leq
\frac{\varepsilon}{|t|}e^{-|n+m|\rho}
\end{eqnarray}

Now we state our abstract KAM theorem, and as a corollary, we get Theorem 1.

\begin{Theorem}
Assume that the Hamiltonian $H=N + {\mathcal B} + P$ satisfies condition $\bf (A1)-(A6)$.
Let $\gamma > 0$ be sufficiently small, then there exists $\varepsilon > 0$ and $\rho > 0$ such that if $\|X_P\|_{D_{\rho}(r,s),\mathcal O} < \varepsilon$, the following holds: There exists a Cantor subset
${\mathcal O}_{\gamma} \subseteq {\mathcal O}$ with $meas({\mathcal
O}\setminus {\mathcal O}_{\gamma}) = O(\gamma^{\varsigma})$
($\varsigma$ is a positive constant) and two maps which are analytic
in $\theta,\varphi$ and $C_W ^4$ in $(\xi,\sigma)$.
$$
\Phi:\T^{b+\tilde{b}} \times {\mathcal O}_{\gamma} \rightarrow D_{\rho}(r,s),\qquad \tilde{\omega}: {\mathcal O}_{\gamma}\rightarrow \R^{b+\tilde{b}}
$$
where $\Phi$ is $\frac{\varepsilon}{\gamma^{16}}$-close to the trivial embedding
$
\Phi_0 : \T^{b+\tilde{b}} \times {\mathcal O} \rightarrow
\T^{b+\tilde{b}} \times \{0,0\} \times \{0,0\} \times \{0,0\}
$
and $\tilde{\omega}$ is $\varepsilon$-close
to the unperturbed frequency $\omega$. Such that
$\forall (\xi,\sigma) \in {\mathcal O}_{\gamma}$ and $(\theta,\varphi) \in \T^{b+\tilde{b}}$, the curve
$
t\rightarrow \Psi \bigg((\theta,\varphi)+ \tilde{\omega}t, (\xi,\sigma) \bigg)
$
is a quasi-periodic solution of the Hamiltonian equation governed by
$H=N+ {\mathcal B} + P$.
\end{Theorem}

\section{\textbf{Normal Form}}
\noindent

\indent Consider the equation set (1.1) in a view of Hamiltonian system and it could be rewritten as
\begin{eqnarray}
\left\{
\begin{array}{l}
-i\partial_t u = \frac{\partial H}{\partial \bar u}  \\
-i\partial_t v = \frac{\partial H}{\partial \bar v}   \\
\end{array}
\right.
\end{eqnarray}
where the function $H$ is a Hamiltonian:
\begin{eqnarray}
H=\frac{1}{2} \bigg(\langle (-\Delta + M_{\xi}) u,u\rangle+\langle (-\Delta + M_{\sigma}) v,v\rangle  \bigg)+
\int_{\T^d}G(|u|^2 , |v|^2) dx
\end{eqnarray}

Expanding $u,v$ into Fourier
series
\begin{eqnarray}
u=\sum_{n\in \Z^d}u_n \phi_n,\quad v=\sum_{n\in \Z^d}v_n \phi_n,\quad \phi_n = \sqrt{\frac{1}{(2\pi)^d}}e^{i\langle n,x\rangle}
\end{eqnarray}
so the Hamiltonian becomes
\begin{eqnarray}
H=\sum_{n \in \Z^d}\lambda_n u_n \bar u_n + \sum_{n \in \Z^d}\tilde{\lambda}_n v_n \bar v_n +
P(\theta,\varphi,I,J,u,\bar u,v,\bar v;\xi,\sigma)
\end{eqnarray}
where $\lambda_n$ is the eigenvalue of $-\Delta+M_{\xi}$ and $\tilde{\lambda}_n$ is the eigenvalue of $-\Delta+M_{\sigma}$, which means
$\lambda_n = |n|^2 + \xi_n$ if $n\in S$ and $\lambda_n = |n|^2$ if $n\in \Z^d _1$;
$\tilde{\lambda}_n = |n|^2 + \sigma_n$ if $n\in \tilde{S}$ and $\tilde{\lambda}_n = |n|^2$ if $n\in \Z^d _2$.

Introducing action-angle variable:
\begin{eqnarray}
&u_n = \sqrt{I_n}e^{i\langle k,\theta\rangle},\bar u_n = \sqrt{I_n}e^{-i\langle k,\theta\rangle}, \quad n\in S \nonumber \\
&v_n = \sqrt{J_n}e^{i\langle k,\varphi\rangle},\bar v_n = \sqrt{J_n}e^{-i\langle k,\varphi\rangle}, \quad n\in \tilde{S}
\end{eqnarray}
The Hamiltonian (3.4) is now turned into
\begin{eqnarray}
H=\langle \omega,I\rangle+\langle \tilde{\omega},J\rangle+
\sum_{n\in \Z^d _1}\Omega_n u_n \bar u_n + \sum_{n\in \Z^d _2}\tilde{\Omega}_n v_n \bar v_n +
P(\theta,\varphi,I,J,u,\bar u,v,\bar v; \xi,\sigma)
\end{eqnarray}

Now let's verify condition $\bf (A1)-(A6)$ for (3.6).

\indent Verifying $\bf (A1)$:
It's easy to see that $\frac{\partial (\omega,\tilde{\omega})}{\partial (\xi,\sigma)} = I_{b+\tilde{b}}$

\indent Verifying $\bf (A2)$:
$\Omega_n=|n|^2, n\in \Z^d _1$ and $\tilde{\Omega}_n=|n|^2, n\in \Z^d _2$ so it's obvious.

\indent Verifying $\bf (A3)$:
For convenience, we only verify the most complicated (2.18). Recall the structure of the Hamiltonian (3.6), now
$a_n =0, b_n =0$, so we only need to concentrate on each element of the diagonal of the matrix
\begin{eqnarray}
\big(
\langle k,\omega \rangle + \langle \tilde{k},\tilde{\omega} \rangle
\big)I + A_n \otimes I_2 - I_2 \otimes A^{\mathrm T} _m
\end{eqnarray}
which means that we only need to verify
\begin{eqnarray}
|\langle k,\omega \rangle + \langle \tilde{k},\tilde{\omega} \rangle + l|\geq \frac{\gamma}{K^{\tau}}
\end{eqnarray}
After excluding a subset with measure $O(\gamma)$ of the parameter set, (3.8) follows and thus condition $\bf (A3)$ is verified.

\indent Verifying $\bf (A4)$:
It's similar with the verification of condition $\bf (A4)$ in $\cite{17}$.

\indent Verifying $\bf (A5)$: It's very similar to that in $\cite{15}$.
Recall the nonlinearity $\int_{\T^d} G(|u|^2, |v|^2) dx$,
expand $G$ into Taylor series in a neighbourhood of the origin and expand $u,v$ into Fourier series, it could be written as a sum of such terms:
$$
u_{a_1}u_{a_2}\cdots u_{a_m}\bar u_{b_1}\bar u_{b_2}\cdots \bar u_{b_m}
v_{c_1}v_{c_2}\cdots v_{c_n}\bar v_{d_1}\bar v_{d_2}\cdots \bar v_{d_n}
$$
with $n,m \geq 1$ and
$$
\sum_{j=1} ^m a_j - \sum_{j=1} ^m b_j + \sum_{j=1} ^n c_j - \sum_{j=1} ^n d_j  = 0
$$
\noindent It belongs to a more general case  :
$$
u^{\lambda}\bar u^{\mu} v^{\tilde{\lambda}} \bar v^{\tilde{\mu}} \quad with \quad
\sum_{n\in \Z^d}(\lambda_n - \mu_n)n + \sum_{n\in \Z^d}(\tilde{\lambda}_n - \tilde{\mu}_n)n = 0 \eqno{(\star)}
$$
we prove that $(\star)$ also satisfies condition $\bf (A5)$. By the definition of $S,\tilde{S}$, we have
\begin{eqnarray}
u^{\lambda}\bar u^{\mu}v^{\tilde{\lambda}}\bar v^{\tilde{\mu}} &=&
\sqrt{I_{i_1}}^{\lambda_{i_1}+\mu_{i_1}}\cdots \sqrt{I_{i_b}}^{\lambda_{i_b}+\mu_{i_b}}
\sqrt{J_{t_1}}^{\tilde{\lambda}_{t_1}+\tilde{\mu}_{t_1}}\cdots \sqrt{J_{t_{\tilde{b}}}}^{\tilde{\lambda}_{t_{\tilde{b}}}+\tilde{\mu}_{t_{\tilde{b}}}}
  \nonumber \\
 &\times&  e^{i\big(\sum_{j=1} ^b (\lambda_{i_j}-\mu_{i_j}) \theta_{i_j} +
             \sum_{j=1} ^{\tilde{b}} (\tilde{\lambda}_{t_j}-\tilde{\mu}_{t_j}) \varphi_{t_j}
\big)}    \nonumber \\
&\times&
u^{\lambda-\sum_{j=1} ^b \lambda_{i_j} e_{i_j}}  \bar{u}^{\mu-\sum_{j=1} ^b \mu_{i_j} e_{i_j}}
v^{\tilde{\lambda}-\sum_{j=1} ^{\tilde{b}} \tilde{\lambda}_{t_j} e_{t_j}}
{\bar v}^{\tilde{\mu}-\sum_{j=1} ^{\tilde{b}} \tilde{\mu}_{t_j} e_{t_j}}
\end{eqnarray}
Just set
$k=(\lambda_{i_1}-\mu_{i_1},\cdots,\lambda_{i_b}-\mu_{i_b}),
\tilde{k}=(\tilde{\lambda}_{t_1}-\tilde{\mu}_{t_1},\cdots,\tilde{\lambda}_{t_{\tilde{b}}}
-\tilde{\mu}_{t_{\tilde{b}}})
$,
$
\alpha=\lambda-\sum_{j=1} ^b \lambda_{i_j} e_{i_j},
\beta=\mu-\sum_{j=1} ^b \mu_{i_j} e_{i_j},
\tilde{\alpha}=\tilde{\lambda}-\sum_{j=1} ^{\tilde{b}} \tilde{\lambda}_{t_j} e_{t_j},
\tilde{\beta}=\tilde{\mu}-\sum_{j=1} ^{\tilde{b}} \tilde{\mu}_{t_j} e_{t_j}
$.
Combining with $(\star)$, it's easy to verify condition $\bf (A5)$.

\indent Verifying $\bf (A6)$:
Now $B=0$. Because $P \in \mathcal A$, we have
\begin{eqnarray}
P=\sum_{kl\alpha\beta, \tilde{k}\tilde{l}\tilde{\alpha}\tilde{\beta}}
   P_{kl\alpha\beta, \tilde{k}\tilde{l}\tilde{\alpha}\tilde{\beta}}(\xi,\sigma)
   e^{i(\langle k,\omega \rangle+\langle \tilde{k},\tilde{\omega} \rangle)}I^l J^{\tilde{l}}
      u^{\alpha}\bar u^{\beta} v^{\tilde{\alpha}} \bar v^{\tilde{\beta}}  \nonumber
\end{eqnarray}
satisfying
\begin{eqnarray}
\sum_{j=1} ^b k_j i_j + \sum_{j=1} ^{\tilde{b}} \tilde{k}_j \tilde{i}_j +
\sum_{n\in \Z^d _1}(\alpha_n - \beta_n)n + \sum_{n\in \Z^d _2}(\tilde{\alpha}_n - \tilde{\beta}_n)n \neq 0
\Rightarrow
P_{kl\alpha\beta, \tilde{k}\tilde{l}\tilde{\alpha}\tilde{\beta}} = 0  \nonumber
\end{eqnarray}
By the fact $(n+tc)-(m+tc)=n-m, (n+tc)+(m-tc)=n+m$, we have
\begin{eqnarray}
&\frac{\partial^2 P}{\partial u_{n} \partial v_{m} } =
\frac{\partial^2 P}{\partial u_{n+tc} \partial v_{m-tc} }
=\lim_{t\rightarrow \infty} \frac{\partial^2 P}{\partial u_{n+tc} \partial v_{m-tc} }  \nonumber \\
&\frac{\partial^2 P}{\partial u_{n} \partial \bar v_{m} } =
\frac{\partial^2 P}{\partial u_{n+tc} \partial \bar v_{m+tc} }
=\lim_{t\rightarrow \infty} \frac{\partial^2 P}{\partial u_{n+tc} \partial \bar v_{m+tc} }  \nonumber \\
&\frac{\partial^2 P}{\partial \bar u_{n} \partial v_{m} } =
\frac{\partial^2 P}{\partial \bar u_{n+tc} \partial v_{m+tc} }
=\lim_{t\rightarrow \infty} \frac{\partial^2 P}{\partial \bar u_{n+tc} \partial v_{m+tc} }  \nonumber \\
&\frac{\partial^2 P}{\partial \bar u_{n} \partial \bar v_{m} } =
\frac{\partial^2 P}{\partial \bar u_{n+tc} \partial \bar v_{m-tc} }
=\lim_{t\rightarrow \infty} \frac{\partial^2 P}{\partial \bar u_{n+tc} \partial \bar v_{m-tc} }  \nonumber
\end{eqnarray}
so (2.23)-(2.27) is verified, and following the same method we could verify (2.19)-(2.22) and (2.28)-(2.31). Thus $\bf (A6)$ is verified.

\section{\textbf{KAM Iteration}}
\noindent

 We prove Theorem 2 by a KAM iteration which involves an infinite sequence of change of variables. Each step of KAM iteration makes the perturbation smaller than the previous step at the cost of excluding a small set of parameters. We have to prove the convergence of the iteration and estimate the measure of the excluded set after infinite KAM steps.

At the $\nu$-th step of the KAM iteration, we consider a Hamiltonian vector field with
$$
H_{\nu} = N_{\nu} + {\mathcal B}_{\nu} + P_{\nu}
(
\theta,\varphi,I,J,u,\bar u,v,\bar v; \xi,\sigma
)
$$
where
\begin{eqnarray}
&&N_{\nu}=\langle \omega(\xi,\sigma),I\rangle+\langle \tilde{\omega}(\xi,\sigma),J\rangle+
\sum_{n\in \Z^d _1}\Omega_n(\xi,\sigma)u_n \bar u_n +
\sum_{n\in \Z^d _2}\tilde{\Omega}_n(\xi,\sigma)v_n \bar v_n \\
&&{\mathcal B}=\sum_{n\in \Z^d _1 \cap \Z^d _2}\bigg(a_n(\xi,\sigma)u_n \bar v_n
+ b_n(\xi,\sigma)\bar u_n v_n \bigg)
\end{eqnarray}
with $B_{\nu} + P_{\nu}$ is defined in $D(r_{\nu}, s_{\nu})\times {\mathcal O}_{\nu -1}$.

We construct a map
$$
\Phi_{\nu}:D(r_{\nu+1},s_{\nu+1})\times {\mathcal O}_{\nu} \rightarrow D(r_{\nu},s_{\nu})\times {\mathcal O}_{\nu-1}
$$
so that the vector field $X_{H_{\nu}\circ \Phi_{\nu}}$ defined on $D(r_{\nu+1},s_{\nu+1})$satisfies
$$
\| X_{P_{\nu+1}} \|_{ D(r_{\nu+1},s_{\nu+1}), {\mathcal O}_{\nu} }=\| X_{H_{\nu}\circ \Phi_{\nu}} -
X_{N_{\nu+1}+{\mathcal B}_{\nu+1}} \|_{ D(r_{\nu+1},s_{\nu+1}), {\mathcal O}_{\nu} } \leq \varepsilon_{\nu} ^{\kappa}, \quad \kappa > 1
$$
and the new Hamiltonian still satisfies condition $\bf (A1)-(A6)$.

\indent To simplify notations, in the following text, the quantities without subscripts refer to quantities at the $\nu$-th step, while the quantities with subscripts $+$ denote the corresponding quantities at the $(\nu+1)$-th step. Let's consider the Hamiltonian defined in $D(r,s)\times \mathcal O$:
\begin{eqnarray}
H&=& N+{\mathcal B} + P \nonumber \\
 &=& \left\langle \omega(\xi,\sigma), I \right\rangle + \langle \tilde{\omega}(\xi,\sigma), J \rangle
 + \sum_{n\in \Z^d_1}\Omega_n(\xi,\sigma) u_n \bar u_n + \sum_{n\in \Z^d_2}\tilde{\Omega}_n(\xi,\sigma) v_n \bar v_n  \nonumber \\
 &+& \sum_{n\in \Z^d _1 \cap \Z^d _2}\bigg(a_n(\xi,\sigma)u_n \bar v_n + b_n(\xi,\sigma)v_n \bar u_n \bigg)
 +P(\theta,\varphi, I,J, u, \bar u,v,\bar v; \xi,\sigma, \varepsilon)
\end{eqnarray}
We assume that for $(\xi,\sigma) \in \mathcal O$, one has:

For any $|k|+|\tilde{k}| \leq K, n \in \Z^d _1, m\in \Z^d _2$, the followings hold
\begin{eqnarray}
| \langle k,\omega \rangle + \langle \tilde{k},\tilde{\omega} \rangle | \geq \frac{\gamma}{K^{\tau}}, \qquad
|k|+|\tilde{k}| \neq 0  \nonumber
\end{eqnarray}
and
\begin{eqnarray}
|\det\left(
\big(\langle k,\omega \rangle + \langle \tilde{k},\tilde{\omega} \rangle \big) I + A_n \right)
|\geq \frac{\gamma}{K^{\tau}}  \nonumber
\end{eqnarray}
and
\begin{eqnarray}
|\det \bigg(
(\langle k,\omega \rangle+\langle \tilde{k},\tilde{\omega} \rangle) I
+ A_n \otimes I_2 - I_2 \otimes A^{\mathrm T} _m
\bigg)| \geq \frac{\gamma}{K^{\tau}},
  \qquad |k|+|\tilde{k}| \neq 0     \nonumber
\end{eqnarray}

Expand $P$ into Fourier-Taylor series
$
P=\sum\limits_{kl\alpha\beta,\tilde{k}\tilde{l}\tilde{\alpha}\tilde{\beta}}
P_{kl\alpha\beta,\tilde{k}\tilde{l}\tilde{\alpha}\tilde{\beta}}
e^{i(\langle k,\theta\rangle+\langle \tilde{k},\varphi\rangle)} I^l J^{\tilde{l}}
u^{\alpha}\bar u^{\beta}v^{\tilde{\alpha}}\bar v^{\tilde{\beta}}
$
and by condition $\bf (A5)$ we get that
\begin{eqnarray}
P_{kl\alpha\beta,\tilde{k}\tilde{l}\tilde{\alpha}\tilde{\beta}}=0 \quad if \quad
\sum_{1\leq j \leq b}k_j i_j + \sum_{n\in \Z^d _1}(\alpha_n - \beta_n)n +
\sum_{1\leq j \leq \tilde{b}}\tilde{k}_j t_j + \sum_{n\in \Z^d _2}(\tilde{\alpha}_n - \tilde{\beta}_n)n
 \neq 0   \nonumber
\end{eqnarray}
which means that when $k=0,\tilde{k}=0$, the terms $u_n\bar u_m, u_n \bar v_m, v_n \bar u_m, v_n\bar v_m$ are absent when $|n|=|m|,n \neq m$.

Now we describe how to construct a subset $\mathcal O_+ \subseteq \mathcal O$ and a change of variables
$\Phi : D_+ \times {\mathcal O}_+ = D(r_+, s_+)\times {\mathcal O}_+ \rightarrow D(r,s)\times {\mathcal O}$ such that the transformed Hamiltonian $H_+ = H\circ \Phi = N_+ + {\mathcal B}_+ + P_+$ satisfies conditions $\bf (A1)-(A6)$ with new parameters $\varepsilon_+, r_+, s_+$ and with $(\xi,\sigma) \in \mathcal O_+$.

\subsection{\textbf{Homological Equation}}
\noindent

Expand $P$ into Fourier-Taylor series
\begin{eqnarray}
P=\sum_{kl\alpha\beta,\tilde{k}\tilde{l}\tilde{\alpha}\tilde{\beta}}
P_{kl\alpha\beta,\tilde{k}\tilde{l}\tilde{\alpha}\tilde{\beta}}
e^{i(\langle k,\theta\rangle+\langle \tilde{k},\varphi\rangle)}I^l J^{\tilde{l}}
u^{\alpha}\bar u^{\beta}v^{\tilde{\alpha}}\bar v^{\tilde{\beta}}
\end{eqnarray}
where $k\in \Z^b, l\in \N^b;\tilde{k}\in \Z^{\tilde{b}},\tilde{l}\in N^{\tilde{b}}$ and the multi-indices $\alpha,\beta;\tilde{\alpha},\tilde{\beta}$
run over the set of all infinite dimensional vectors $\alpha =
(\cdots ,\alpha_n,\cdots)_{n\in \Z^d _1},\beta = (\cdots
,\beta_n,\cdots)_{n\in \Z^d _1};
\tilde{\alpha} = (\cdots ,\tilde{\alpha}_n,\cdots)_{n\in \Z^d _2},
\tilde{\beta} = (\cdots,\tilde{\beta}_n,\cdots)_{n\in \Z^d _2}
$
with finitely many nonzero
components of positive integers. And by $\bf (A5)$ we get that
\begin{eqnarray}
P_{kl\alpha\beta,\tilde{k}\tilde{l}\tilde{\alpha}\tilde{\beta}}=0 \quad if \quad
\sum_{1\leq j \leq b}k_j i_j + \sum_{1\leq j \leq \tilde{b}}\tilde{k}_j t_j +
\sum_{n\in \Z^d _1}(\alpha_n - \beta_n)n + \sum_{n\in \Z^d _2}(\tilde{\alpha}_n - \tilde{\beta}_n)n
 \neq 0
\end{eqnarray}

Consider its quadratic truncation $R$:
\begin{eqnarray}
R(\theta,\varphi, I,J, u,\bar u,v,\bar v)
= R_0 + R_1 + R_2
\end{eqnarray}
where
\begin{eqnarray}
R^0 = \sum_{|k|+|\tilde{k}| \leq K,|l|+|\tilde{l}|\leq 1}
 \left(P^{k100,\tilde{k}000} _l I^l +
 P^{k000,\tilde{k}100} _{\tilde{l}} J^{\tilde{l}}   \right)
 e^{i(\langle k,\theta \rangle+ \langle \tilde{k},\varphi\rangle)}
\end{eqnarray}

\begin{eqnarray}
R^1=\sum_{|k|+|\tilde{k}| \leq K;n}
\left(
P^{k10,\tilde{k}00} _n u_n + P^{k01,\tilde{k}00} _n \bar u_n +
P^{k00,\tilde{k}10} _n v_n + P^{k00,\tilde{k}01} _n \bar v_n
\right)
e^{i(\langle k,\theta \rangle+ \langle \tilde{k},\varphi\rangle)}
\end{eqnarray}
and
\begin{eqnarray}
R^2=R^2 _{uu} + R^2 _{uv} + R^2 _{vv}
\end{eqnarray}
where
\begin{eqnarray}
R^2 _{uu}=\sum_{|k|+|\tilde{k}| \leq K;nm}\left(P^{k20,\tilde{k}00} _{nm}u_n u_m +
          P^{k11,\tilde{k}00} _{nm}u_n \bar u_m +
         P^{k02,\tilde{k}00} _{nm}\bar u_n \bar u_m \right)
         e^{i(\langle k,\theta \rangle+\langle \tilde{k},\varphi \rangle)}
\end{eqnarray}

\begin{eqnarray}
R^2 _{uv} &=& \sum_{|k|+|\tilde{k}| \leq K;nm}
          \bigg(P^{k10,\tilde{k}10} _{nm}u_n v_m +
          P^{k10,\tilde{k}01} _{nm}u_n \bar v_m
          + P^{k01,\tilde{k}10} _{nm}\bar u_n v_m   \nonumber \\
          &+& P^{k01,\tilde{k}01} _{nm}\bar u_n \bar v_m  \bigg)
           e^{i(\langle k,\theta \rangle+\langle \tilde{k},\varphi \rangle)}
\end{eqnarray}

\begin{eqnarray}
R^2 _{vv}=\sum_{|k|+|\tilde{k}| \leq K;nm}\left(P^{k00,\tilde{k}20} _{nm}v_n v_m +
          P^{k00,\tilde{k}11} _{nm}v_n \bar v_m +
         P^{k00,\tilde{k}02} _{nm}\bar v_n \bar v_m \right)
         e^{i(\langle k,\theta \rangle+\langle \tilde{k},\varphi \rangle)}
\end{eqnarray}

We explain the coefficients in (4.7)-(4.12) as below (here $e_n$ denotes the vector with the $n$ th component being $1$ and the other components being zero ):

$P^{k100,\tilde{k}000} _l = P_{kl00,\tilde{k}000},P^{k000,\tilde{k}100} _{\tilde{l}} = P_{k000,\tilde{k}\tilde{l}00} $;

$P^{k10,\tilde{k}00}_n = P_{k0\alpha\beta,\tilde{k}0\tilde{\alpha}\tilde{\beta}}$ with
$\alpha = e_n,\beta = 0;\tilde{\alpha}=0,\tilde{\beta} = 0$;

$P^{k01,\tilde{k}00}_n = P_{k0\alpha\beta,\tilde{k}0\tilde{\alpha}\tilde{\beta}}$ with
$\alpha = 0,\beta = e_n;\tilde{\alpha}=0,\tilde{\beta} = 0$;

the definitions of $P^{k00,\tilde{k}10}_n,P^{k00,\tilde{k}01}_n$ are similar.

$P^{k20,\tilde{k}00}_{nm} = P_{k0\alpha\beta,\tilde{k}0\tilde{\alpha}\tilde{\beta}}$ with
$\alpha=e_n + e_m,\beta=0,\tilde{\alpha}=0,\tilde{\beta}=0$;

$P^{k11,\tilde{k}00}_{nm} = P_{k0\alpha\beta,\tilde{k}0\tilde{\alpha}\tilde{\beta}}$ with
$\alpha=e_n,\beta=e_m,\tilde{\alpha}=0,\tilde{\beta}=0$;

$P^{k02,\tilde{k}00}_{nm} = P_{k0\alpha\beta,\tilde{k}0\tilde{\alpha}\tilde{\beta}}$ with
$\alpha=0,\beta=e_n + e_m,\tilde{\alpha}=0,\tilde{\beta}=0$;

the definitions of
$P^{k00,\tilde{k}20}_{nm},P^{k00,\tilde{k}11}_{nm},P^{k00,\tilde{k}02}_{nm}$ are similar.

$P^{k10,\tilde{k}10}_{nm} = P_{k0\alpha\beta,\tilde{k}0\tilde{\alpha}\tilde{\beta}}$ with
$\alpha=e_n,\beta=0,\tilde{\alpha}=e_m,\tilde{\beta}=0$;

$P^{k10,\tilde{k}01}_{nm} = P_{k0\alpha\beta,\tilde{k}0\tilde{\alpha}\tilde{\beta}}$ with
$\alpha=e_n,\beta=0,\tilde{\alpha}=0,\tilde{\beta}=e_m$;

the definitions of
$P^{k01,\tilde{k}10}_{nm},P^{k01,\tilde{k}01}_{nm} $ are similar.

 $\\$
\indent Rewrite $H$ as $H= N+ {\mathcal B} + R + (P - R)$. Due to the choice of $s_+ \ll s$ (defined in section 5) and the definition of the norm, it follows  immediately
\begin{eqnarray}
\| X_R \|_{D(r,s),\mathcal O}\leq \| X_P \|_{D(r,s),\mathcal O}\leq \varepsilon
\end{eqnarray}
and in $D(r,s_+)$
\begin{eqnarray}
\| X_{P-R} \|_{D(r,s_+)}\leq \varepsilon_+
\end{eqnarray}

In the following, we will construct a Hamiltonian function $F$ satisfying $\bf (A5)$ and with the same form of $R$ defined in $D_+ = D(r_+ , s_+)$ such that the time one map $X^1 _F$ of the Hamiltonian vector field $X_F$ defines a map from $D_+$ to $D$ and puts $H$ into $H_+$. Precisely, one has
\begin{eqnarray}
H\circ X^1 _F
&=& (N+{\mathcal B}+ R)\circ X^1 _F + (P-R)\circ X^1 _F \nonumber \\
&=& N+ {\mathcal B} \nonumber \\
&+& \{N+  {\mathcal B} , F \} + R  \\
&+& \int_0 ^1 (1-t)\{ \{N+{\mathcal B}, F \}, F \}\circ X^t _F dt \\
&+& \int_0 ^1 \{R, F\}\circ X^t _F dt + (P-R)\circ X^1 _F
\end{eqnarray}
So we get the linearized homological equation:
\begin{eqnarray}
\{N+ {\mathcal B}, F \}+ R = {\hat N} + \hat{\mathcal B}
\end{eqnarray}
where
\begin{eqnarray}
&&\hat N = P_{0000,0000} + \langle \hat{\omega},I\rangle + \langle \hat{\tilde{\omega}},J\rangle
+ \sum_{n\in \Z^d _1}P^{011,000}_{nn} u_n \bar u_n
+ \sum_{n\in \Z^d _2}P^{000,011}_{nn} v_n \bar v_n  \\
&& \hat{\omega} = (P^{0100,0000} _l)_{|l|=1}, \quad \hat{\tilde{\omega}}=(P^{0000,0100} _{\tilde{l}})_{|\tilde{l}|=1} \\
&& \hat{\mathcal B} = \sum_{n\in \Z^d _1 \cap \Z^d _2}
(P^{010,001} _{nn} u_n \bar v_n + P^{001,010} _{nn} \bar u_n v_n)
\end{eqnarray}
We define
$
N_+ = N+\hat N, {\mathcal B}_+ = {\mathcal B} + \hat{\mathcal B},
$
and
\begin{eqnarray}
P_+ = \int_0 ^1 (1-t)\{ \{N+ {\mathcal B}, F \}, F \}\circ X^t _F dt
+ \int_0 ^1 \{R, F\}\circ X^t _F dt + (P-R)\circ X^1 _F
\end{eqnarray}
We construct the Hamiltonian function $F$ as below, with the same structure of $R$:
\begin{eqnarray}
F(\theta,\varphi, I,J,u ,\bar z , v,\bar v)=F^0 + F^1 + F^2
\end{eqnarray}
with
\begin{eqnarray}
F^0 = \sum_{|k|+|\tilde{k}|\leq K,|l|+|\tilde{l}| \leq 1}\left(F^{k100,\tilde{k}000} I^l + F^{k000,\tilde{k}100} J^{\tilde{l}} \right)
e^{i(\langle k,\theta \rangle + \langle \tilde{k},\varphi \rangle )}
\end{eqnarray}

\begin{eqnarray}
F^1 = \sum_{|k|+|\tilde{k}|\leq K,n}\left(
F^{k10,\tilde{k}00} _n u_n + F^{k01,\tilde{k}00} _n \bar u_n +
F^{k00,\tilde{k}10} _n v_n + F^{k00,\tilde{k}01} _n \bar v_n
\right)
e^{i(\langle k,\theta \rangle+\langle \tilde{k},\varphi\rangle)}
\end{eqnarray}

and
\begin{eqnarray}
F^2 = F^2 _{uu} + F^2 _{uv} + F^2 _{vv}
\end{eqnarray}
where
\begin{eqnarray}
F^2 _{uu}=\sum_{|k|+|\tilde{k}|\leq K,nm}\left(
F^{k20,\tilde{k}00} _{nm}u_n u_m +
F^{k11,\tilde{k}00} _{nm}u_n \bar u_m +
F^{k02,\tilde{k}00} _{nm}\bar u_n \bar u_m
\right)
e^{i(\langle k,\theta \rangle+\langle \tilde{k},\varphi\rangle)}
\end{eqnarray}

\begin{eqnarray}
 F^2 _{uv}&=&\sum_{|k|+|\tilde{k}|\leq K,nm}
\bigg(
F^{k10,\tilde{k}10} _{nm}u_n v_m +
F^{k10,\tilde{k}01} _{nm}u_n \bar v_m +
F^{k01,\tilde{k}10} _{nm}\bar u_n v_m  \nonumber \\
&+& F^{k01,\tilde{k}01} _{nm}\bar u_n \bar v_m
\bigg)
e^{i(\langle k,\theta \rangle+\langle \tilde{k},\varphi\rangle)}
\end{eqnarray}

\begin{eqnarray}
F^2 _{vv}=\sum_{|k|+|\tilde{k}|\leq K,nm}\bigg(
F^{k00,\tilde{k}20} _{nm}v_n v_m +
F^{k00,\tilde{k}11} _{nm}v_n \bar v_m +
F^{k00,\tilde{k}02} _{nm}\bar v_n \bar v_m
\bigg)
e^{i(\langle k,\theta \rangle+\langle \tilde{k},\varphi\rangle)}
\end{eqnarray}

Now (4.18) is turned into
\begin{eqnarray}
&\{N, F^0\}+ R^0 - P_{0000,0000} - \langle \hat{\omega},I \rangle - \langle \hat{\tilde{\omega}},J \rangle = 0  \\
&\{N+ {\mathcal B}, F^1\}+ R^1 = 0
\end{eqnarray}
and the most complicated
\begin{eqnarray}
\{N + {\mathcal B}, F^2\}+R^2
 &=& \sum_{n\in \Z^d _1}P^{011,000}_{nn}u_n \bar u_n
+  \sum_{n\in \Z^d _2}P^{000,011}_{nn}u_n \bar u_n   \nonumber \\
&+&
\sum_{n\in \Z^d _1 \cap \Z^d _2}
( P^{010,001} _{nn} u_n \bar v_n + P^{001,010} _{nn} \bar u_n v_n )
\end{eqnarray}
Now we solve equations (4.30)-(4.32) one by one:

$\bf Solving (4.30)$:
By the expansion (4.24), (4.30) is turned into
\begin{eqnarray}
&i(\langle k,\omega\rangle+\langle \tilde{k},\tilde{\omega}\rangle) F^{k100,\tilde{k}000}
=P^{k100,\tilde{k}000}, \qquad |k|+|\tilde{k}|\neq 0  \\
&i(\langle k,\omega\rangle+\langle \tilde{k},\tilde{\omega}\rangle) F^{k000,\tilde{k}100}
=P^{k000,\tilde{k}100}, \qquad |k|+|\tilde{k}| \neq 0
\end{eqnarray}
According to assumption (2.15) in condition $\bf (A3)$, one has
\begin{eqnarray}
&|F^{k100,\tilde{k}000}|_{\mathcal O} \leq \gamma^{-16} K^{16\tau+16} |F^{k100,\tilde{k}000}|_{\mathcal O},
\quad 0< |k|+|\tilde{k}|\leq K    \\
&|F^{k000,\tilde{k}100}|_{\mathcal O} \leq \gamma^{-16} K^{16\tau+16} |F^{k000,\tilde{k}100}|_{\mathcal O},
\quad 0< |k|+|\tilde{k}|\leq K
\end{eqnarray}

$\bf Solving (4.31)$: For convenience, we only describe the homological equation (with the estimate of the solution) related to the elimination of term $u_n,  v_n, n\in \Z^d _1 \cap \Z^d _2$, and the corresponding equation related to $\bar u_n, \bar v_n$ with the estimate of its solutions follow the same way.
By the expansion (4.25), (4.31) is turned into
\begin{eqnarray}
& \big( \langle k,\omega \rangle+\langle \tilde{k},\tilde{\omega} \rangle + \Omega_n \big)F^{k10,\tilde{k}00} _n
 +a_n F^{k00,\tilde{k}10} _n = P^{k10,\tilde{k}00} _n   \\
& \big( \langle k,\omega \rangle+\langle \tilde{k},\tilde{\omega} \rangle + \tilde{\Omega}_n \big)F^{k00,\tilde{k}10} _n
 +b_n F^{k00,\tilde{k}10} _n = P^{k00,\tilde{k}10} _n
\end{eqnarray}
The coefficient matrix is just $(\langle k,\omega \rangle + \langle \tilde{k},\tilde{\omega} \rangle)I + A_n$,
so according to assumption (2.16) in condition $\bf (A3)$, one has the estimate
\begin{eqnarray}
|F^{k10,\tilde{k}00} _n|_{\mathcal O}, |F^{k00,\tilde{k}10} _n|_{\mathcal O} <
\gamma^{-16} K^{16\tau+16} \varepsilon e^{-(|k|+|\tilde{k}|)r} e^{-|n|\rho}  \qquad |k|+|\tilde{k}| \leq K
\end{eqnarray}

For the case when $n\in \Z^d _1 \setminus \Z^d _2$, we eliminate the term $u_n$, and the corresponding equation in (4.31) is:
\begin{eqnarray}
(\langle k,\omega\rangle+\langle \tilde{k},\tilde{\omega}\rangle + \Omega_n)F^{k10,\tilde{k}00} _n = P^{k10,\tilde{k}00} _n
\end{eqnarray}
the coefficient matrix is still $(\langle k,\omega\rangle+\langle \tilde{k},\tilde{\omega}\rangle)I+A_n$ by the definition of $A_n$. The case $n\in \Z^d _2 \setminus \Z^d _1$ is similar. So (4.39) still holds.

$\bf Solving (4.32)$: For convenience, we only describe the most complicated equation (with the estimate of its solutions) related to the elimination of terms:
$
u_n \bar u_m, u_n \bar v_m, v_n \bar u_m, v_n \bar v_m; n,m\in \Z^d _1 \cap \Z^d _2
$
with $|n|\neq |m|$. The corresponding homological equation is
\begin{eqnarray}
\left((\langle k,\omega \rangle+\langle \tilde{k},\tilde{\omega} \rangle)I+M\right)X=Y
\end{eqnarray}
where
\begin{eqnarray}
&X=(F^{k11,\tilde{k}00} _{nm}, F^{k10,\tilde{k}01} _{nm}, F^{k01,\tilde{k}10} _{mn}, F^{k00,\tilde{k}11} _{nm}
)^{\mathrm T}  \\
&Y=(P^{k11,\tilde{k}00} _{nm}, P^{k10,\tilde{k}01} _{nm}, P^{k01,\tilde{k}10} _{mn}, P^{k00,\tilde{k}11} _{nm}
)^{\mathrm T}  \\
\end{eqnarray}
and the $4 \times 4$ matrix $M$ is
\begin{eqnarray}
\left(
\begin{array}{cccc}
-\Omega_n +\Omega_m & b_m  & -a_n & 0          \\
a_m & -\Omega_n + \tilde{\Omega}_m & 0 & -a_n \\
-b_n & 0 & -\tilde{\Omega}_n + \Omega_m & b_m \\
0 & -b_n & a_m & -\tilde{\Omega}_n + \tilde{\Omega}_m
\end{array}
\right)
\end{eqnarray}
so the coefficient matrix of (4.41) is just
$
(\langle k,\omega \rangle+\langle \tilde{k},\tilde{\omega} \rangle)I +
A_n \otimes I_2 - I_2 \otimes A^{\mathrm T} _m
$,
according to assumption (2.17) in condition $\bf (A3)$, one has
\begin{eqnarray}
|F^{k11,\tilde{k}00} _{nm}|_{\mathcal O}, |F^{k10,\tilde{k}01} _{nm}|_{\mathcal O},
|F^{k01,\tilde{k}10} _{mn}|_{\mathcal O}, |F^{k00,\tilde{k}11} _{nm}|_{\mathcal O}
\leq \gamma^{-16}K^{16\tau+16} \varepsilon e^{-(|k|+|\tilde{k}|)r} e^{-|n-m|\rho}
\end{eqnarray}
For equation concerning elimination of terms
$ u_n u_m, u_n v_m, v_n u_m, v_n v_m; n,m\in \Z^d _1 \cap \Z^d _2 $, and the equation concerning elimination of terms
$ \bar u_n \bar u_m, \bar u_n \bar v_m, \bar v_n \bar u_m, \bar v_n \bar v_m; n,m\in \Z^d _1 \cap \Z^d _2 $, similar estimates follow by making use of (2.17).

When $n\in \Z^d _1 \setminus \Z^d _2, m \in \Z^d _1 \cap \Z^d _2$, we have the linear equation concerning the elimination of $u_n \bar v_m, u_n \bar u_m$:
\begin{eqnarray}
& (\langle k,\omega\rangle+\langle \tilde{k},\tilde{\omega}\rangle - \Omega_n + \tilde{\Omega}_m)F^{k10,\tilde{k}01} _{nm} + a_m F^{k11,\tilde{k}00} _{nm} = P^{k10,\tilde{k}01} _{nm}  \\
& (\langle k,\omega\rangle+\langle \tilde{k},\tilde{\omega}\rangle - \Omega_n + \Omega_m)F^{k11,\tilde{k}00} _{nm} + b_m F^{k10,\tilde{k}01} _{nm} = P^{k11,\tilde{k}00} _{nm}
\end{eqnarray}
the coefficient matrix is still in the form $(\langle k,\omega\rangle+\langle \tilde{k},\tilde{\omega}\rangle)I + A_n \otimes I_2 - I_2 \otimes {A_m}^{\mathrm T}$. So the same estimates still hold. Similarly we could work on the case $n\in \Z^d _1 \cap \Z^d _2, m\in \Z^d _2 \setminus \Z^d _1$ and $n\in \Z^d _1 \setminus \Z^d _2, m\in \Z^d _2 \setminus \Z^d _1$.

\subsection{\textbf{Estimate of transformation of coordinates}}
\noindent

Now we give the estimate of $X_F$ and $\phi^1 _F$.
\begin{Lemma}
Let $D_i = D(r_+ + \frac{i}{4}(r-r_+),\frac{i}{4}s), 0<i\leq 4$, then one has
\begin{eqnarray}
\|X_F\|_{D_3,\mathcal O} \leq c\gamma^{-16}K^{16\tau+16}\varepsilon
\end{eqnarray}
\end{Lemma}

\begin{Lemma}
Let $\eta=\varepsilon^{\frac{1}{3}}, D_{i\eta} = D(r_+ + \frac{i}{4}(r-r_+),\frac{i}{4}\eta s), 0<i\leq 4$. If
$\varepsilon \ll \gamma^{24} K^{-24\tau-24}$, one has
\begin{eqnarray}
\phi^t _F : D_{2\eta}\rightarrow D_{3\eta}, \quad -1 \leq t \leq 1
\end{eqnarray}
and
\begin{eqnarray}
\|D\phi^t _F - Id\|_{D_{1\eta}} \leq c\gamma^{-16} K^{16\tau+16} \varepsilon
\end{eqnarray}
\end{Lemma}

The proof of Lemma 4.1 and Lemma 4.2 is similar to Lemma 4.2 and Lemma 4.3 in $\cite{17}$ so we omit it.
\subsection{\textbf{Estimate of new perturbations}}
\noindent

Recall the definition of new perturbation
\begin{eqnarray}
&P_+ = \int_0 ^1 \{R(t),F\}\circ \phi^t _F dx + (P-R)\circ \phi^1 _F   \nonumber
\end{eqnarray}
with $R(t)= (1-t)(\{N+{\mathcal B},F\}+R)+tR $. So
\begin{eqnarray}
X_{P_+}=\int_0 ^1 (\phi^t _F)^{\ast} X_{\{R(t),F\}} dt + (\phi^1 _F)^{\ast} X_{P-R}  \nonumber
\end{eqnarray}
\indent By Lemma 4.1, one has
\begin{eqnarray}
\|D\phi^t _F - Id\|_{D_{1\eta}} \leq c\gamma^{-16}K^{16\tau+16}\varepsilon, \quad -1 \leq t \leq 1  \nonumber
\end{eqnarray}
so we get $\|D\phi^t _F\|_{D_{1\eta}} \leq 2, -1\leq t \leq 1$, and
$
\|X_{\{R(t),F\}}\|_{D_{2\eta}} \leq c\gamma^{-16}K^{16\tau+16}\eta^{-2} \varepsilon^2
$. Combining with
$
\|X_{P-R}\|_{D_{2\eta}} \leq c\eta \varepsilon
$, we have the estimate
$$
\|X_{P_+}\|_{D(r_+,s_+)} \leq
c\eta \varepsilon + c\gamma^{-16}K^{16\tau+16}\eta^{-2} \varepsilon^2 < \varepsilon_+
$$

\subsection{\textbf{Estimate of new normal form}}
\noindent

Due to the special form of $P$ defined in $\bf (A5)$, the terms
$
u_n \bar u_m, u_n \bar v_m, v_n \bar u_m, v_n \bar v_m
$ with $|n|=|m|, n\neq m$ are absent, which means that our normal form has a simpler form (We omit the constant term in the normal form part):
\begin{eqnarray}
 N_+ &=& N+\langle \hat{\omega}, I \rangle + \langle \hat{\tilde{\omega}}, J \rangle
     + \sum_{n\in \Z^d _1} P^{011,000} _{nn}u_n \bar u_n
     + \sum_{n\in \Z^d _2} P^{000,011} _{nn}v_n \bar v_n  \nonumber \\
  &+& \sum_{n\in \Z^d _1 \cap \Z^d _2} (P^{010,001} _{nn}u_n \bar v_n + P^{001,010} _{nn}\bar u_n v_n )
 \nonumber \\
 &=& \langle \omega_+,I \rangle + \langle \tilde{\omega}_+,J \rangle +
    \sum_{n\in \Z^d _1}\Omega_n ^+ u_n \bar u_n +
    \sum_{n\in \Z^d _2}\tilde{\Omega}_n ^+ v_n \bar v_n +
    \sum_{n\in \Z^d _1 \cap \Z^d _2} (a_n ^+ u_n \bar v_n + b_n ^+ \bar u_n v_n )  \nonumber
\end{eqnarray}
where
\begin{eqnarray}
&& \omega_+ = \omega + P_{0100,0000}, \quad \tilde{\omega}_+ = \tilde{\omega} + P_{0000,0100};  \nonumber \\
&&   \Omega_n ^+ = \Omega_n + P^{011,000} _{nn}, n\in \Z^d _1; \quad
    \tilde{\Omega}_n ^+ = \tilde{\Omega}_n + P^{011,000} _{nn}, n\in \Z^d _2; \nonumber \\
&& a_n ^+ = a_n + P^{010,001} _{nn}, b_n ^+ = b_n + P^{001,010} _{nn}, \quad n\in \Z^d _1 \cap \Z^d _2 \nonumber
\end{eqnarray}

So with the help of the regularity of $X_P$ and Cauchy estimates, we have
$$
|\omega_+ - \omega|, |\tilde{\omega}_+ - \tilde{\omega}|,
|\Omega_+ - \Omega|, |\tilde{\Omega}_+ - \tilde{\Omega}|,
|a_n ^+ - a_n|, |b_n ^+ - b_n|<\varepsilon
$$
It follows that for $|k|+|\tilde{k}| \leq K$, we have
\begin{eqnarray}
&|\langle k,\omega_+\rangle + \langle \tilde{k},\tilde{\omega}_+\rangle| \geq
|\langle k,\omega\rangle + \langle \tilde{k},\tilde{\omega}\rangle| - \varepsilon K
> \frac{\gamma}{K_+ ^{\tau}}  \nonumber \\
&|\det \bigg((\langle k,\omega_+\rangle + \langle \tilde{k},\tilde{\omega}_+\rangle)I+A_n ^+ \bigg)| \geq
|\det \bigg((\langle k,\omega \rangle + \langle \tilde{k},\tilde{\omega} \rangle)I+A_n  \bigg)| -
\varepsilon K > \frac{\gamma}{K_+ ^{\tau}}  \nonumber
\end{eqnarray}
and
\begin{eqnarray}
&|\det \bigg((\langle k,\omega_+\rangle + \langle \tilde{k},\tilde{\omega}_+\rangle)I
+A_n ^+ \otimes I_2 - I_2 \otimes {A_m ^+}^{\mathrm T}  \bigg)|  \nonumber \\
  & \geq |\det \bigg((\langle k,\omega \rangle + \langle \tilde{k},\tilde{\omega} \rangle)I
   +A_n \otimes I_2 - I_2 \otimes A_m ^{\mathrm T} \bigg)| -
   \varepsilon K > \frac{\gamma}{K_+ ^{\tau}}  \nonumber
\end{eqnarray}
if $|k|+|\tilde{k}|+||n|-|m|| \neq 0$.
Thus in the next KAM iterative step, small denominator conditions are automatically satisfied for
$|k|+|\tilde{k}| \leq K$.

\subsection{\textbf{Verifying of condition $\bf (A5)$ and $\bf (A6)$ after transformation }}
\noindent

Let's verify condition $\bf (A5)$ now. Certainly it could be done by giving a tedious proof which could be modified from Lemma 4.4 in $\cite{15}$, while we prefer to give a much simpler proof here, by making use of the property of Poisson Bracket of two monomials. We need to prove that $P_+ \in \mathcal A$, recall its definition (4.22) and it could be rewritten as:
\begin{eqnarray}
P_+ &=& P-R +\{P,F\}+\frac{1}{2!}\{\{N+{\mathcal B},F\},F\}+ \frac{1}{2!}\{\{P,F\},F\}  \nonumber \\
&+& \cdots +\frac{1}{n!}\{\cdots\{N+{\mathcal B},F\}\cdots,F\}+\cdots +\frac{1}{n!}\{\cdots\{P,F\}\cdots,F\} + \cdots
\end{eqnarray}
Since $P \in \mathcal A$, we have $R,P-R \in \mathcal A$. By the definition of $F$ one has $F \in \mathcal A$. Now we only need to prove that the second line of (4.53) is also in $\mathcal A$. We know $N,{\mathcal B},F,P \in {\mathcal A}$, so we only need to prove the following Lemma:
\begin{Lemma}
If $G(\theta,\varphi,I,J,u,\bar u,v,\bar v),H(\theta,\varphi,I,J,u,\bar u,v,\bar v) \in \mathcal A$, then
$B(\theta,\varphi,I,J,u,\bar u,v,\bar v) = \{G,H\} \in \mathcal A$.
\end{Lemma}
\noindent $\bf Proof:$
$G,H$ are sums of a series of monomials with the form:
\begin{eqnarray}
G_{kl\alpha\beta,\tilde{k}\tilde{l}\tilde{\alpha}\tilde{\beta}}
e^{i(\langle k,\theta\rangle+\langle \tilde{k},\varphi\rangle)} I^l J^{\tilde{l}}
u^{\alpha}\bar u^{\beta} v^{\tilde{\alpha}} \bar v^{\tilde{\beta}}, \quad
H_{pq\gamma\delta,\tilde{p}\tilde{q}\tilde{\gamma}\tilde{\delta}}
e^{i(\langle p,\theta\rangle+\langle \tilde{p},\varphi\rangle)} I^q J^{\tilde{q}}
u^{\gamma}\bar u^{\delta} v^{\tilde{\gamma}} \bar v^{\tilde{\delta}}
\end{eqnarray}
satisfying
\begin{eqnarray}
&&\sum_{j=1} ^b k_j i_j + \sum_{j=1} ^{\tilde{b}} \tilde{k}_j t_j +\sum_{n\in \Z^d _1}(\alpha_n-\beta_n)n + \sum_{n\in \Z^d _2}(\tilde{\alpha}_n-\tilde{\beta}_n)n = 0  \\
&&\sum_{j=1} ^b p_j i_j + \sum_{j=1} ^{\tilde{b}} \tilde{p}_j t_j +
\sum_{n\in \Z^d _2}(\gamma_n-\delta_n)n + \sum_{n\in \Z^d _2}(\tilde{\gamma}_n-\tilde{\delta}_n)n = 0
\end{eqnarray}
Recall the definition of Poisson Bracket:
\begin{eqnarray}
\{G,H\}&=& \langle \frac{\partial G}{\partial \theta}, \frac{\partial H}{\partial I}\rangle -
           \langle \frac{\partial G}{\partial I}, \frac{\partial H}{\partial \theta}\rangle +
           \langle \frac{\partial G}{\partial \varphi}, \frac{\partial H}{\partial J}\rangle -
           \langle \frac{\partial G}{\partial J}, \frac{\partial H}{\partial \varphi}\rangle \nonumber \\
     &+&  i\sum_{n\in \Z^d _1}\big(
      \frac{\partial G}{\partial u_n}\frac{\partial H}{\partial \bar u_n} -
      \frac{\partial G}{\partial \bar u_n}\frac{\partial H}{\partial u_n}
      \big)
     + i\sum_{n\in \Z^d _2}\big(
      \frac{\partial G}{\partial v_n}\frac{\partial H}{\partial \bar v_n} -
      \frac{\partial G}{\partial \bar v_n}\frac{\partial H}{\partial v_n}
      \big)
\end{eqnarray}

 Let's consider the Poisson Bracket of the two terms in (4.54) and for convenience, we just omit the coefficients
( assume the coefficients to be 1). Denote the first one of (4.54) by $A$ and the second one $B$. Obviously $AB \in {\mathcal A}$. We have
$
\frac{\partial A}{{\partial \theta}_{i_j}} \frac{\partial B}{{\partial I}_{i_j}}=
ik_{i_j} q_{i_j}I_{i_j} ^{-1}AB
$ if $q_{i_j} > 0$, and be $0$ if $q_{i_j} = 0$. So we conclude that
$
\langle \frac{\partial G}{\partial \theta},\frac{\partial H}{\partial I}\rangle \in \mathcal A
$
and similarly we have
$
\langle \frac{\partial G}{\partial I},\frac{\partial H}{\partial \theta}\rangle,
\langle \frac{\partial G}{\partial \varphi},\frac{\partial H}{\partial J}\rangle,
\langle \frac{\partial G}{\partial J},\frac{\partial H}{\partial \varphi}\rangle
$
are also in $\mathcal A$. For the remaining terms,
$
\frac{\partial G}{\partial u_n}\frac{\partial H}{\partial \bar u_n} =
\alpha_n \delta_n (u_n \bar u_n)^{-1}AB
$
if $\alpha_n,\delta_n >0$ and be $0$ otherwise. So we conclude that
$
\frac{\partial G}{\partial u_n}\frac{\partial H}{\partial \bar u_n} \in \mathcal A
$
and similarly,
$
\frac{\partial G}{\partial \bar u_n}\frac{\partial H}{\partial u_n},
\frac{\partial G}{\partial v_n}\frac{\partial H}{\partial \bar v_n},
\frac{\partial G}{\partial \bar v_n}\frac{\partial H}{\partial v_n} \in \mathcal A
$. To sum up, we have $\{G,H\} \in \mathcal A$. \qquad $\Box$

$\\$
\indent By Lemma 4.3, the conclusion $P_+ \in \mathcal A$ follows.

$\\$
\indent Now let's turn to verify that $P_+$ satisfies condition $\bf (A6)$ with $K_+, \varepsilon_+, \rho_+$ in place of $K, \varepsilon, \rho$. Recall that
\begin{eqnarray}
P_+ &=& P-R +\{P,F\}+\frac{1}{2!}\{\{N+{\mathcal B},F\},F\}+ \frac{1}{2!}\{\{P,F\},F\}  \nonumber \\
&+& \cdots +\frac{1}{n!}\{\cdots\{N+{\mathcal B},F\}\cdots,F\}+\cdots +\frac{1}{n!}\{\cdots\{P,F\}\cdots,F\} + \cdots \nonumber
\end{eqnarray}
then for a fixed $c \in \Z^d \setminus \{0\}$, and $|n-m|\geq K$ with
$ K\geq \frac{1}{\rho-\rho_+} \ln (\frac{\varepsilon}{\varepsilon_+}) $, one has
\begin{eqnarray}
\|  \frac{\partial^2 (P-R)}{\partial u_{n+tc} \partial \bar v_{m+tc}} -
\lim_{t\rightarrow \infty} \frac{\partial^2 (P-R)}{\partial u_{n+tc} \partial \bar v_{m+tc}}  \|  \leq
\frac{\varepsilon}{|t|} e^{-|n-m|\rho} \leq \frac{\varepsilon_+}{|t|} e^{-|n-m|\rho_+}   \nonumber
\end{eqnarray}
So we get that $P-R$ satisfies $\bf (A6)$ with $K_+,\varepsilon_+,\rho+$ in place of $K,\varepsilon,\rho$. To verify that other terms also satisfy $\bf (A6)$, we only need the following two lemmas.
\begin{Lemma}
$F$ satisfies $\bf (A6)$ with $\varepsilon^{\frac{2}{3}}$ in place of $\varepsilon$.
\end{Lemma}

\noindent$\bf Proof:$ This proof is just modified from Lemma 4.3 in $\cite{12}$ and actually they are very similar.
Given $n,m \in \Z^d _1 \cap \Z^d _2 $, we let
\begin{eqnarray}
\left\{
\begin{array}{l}
Q_{k\tilde{k},nm}=\left(
F^{k11,\tilde{k}00} _{nm}, F^{k10,\tilde{k}01} _{nm}, F^{k01,\tilde{k}10} _{mn}, F^{k00,\tilde{k}11} _{nm}
\right)^{\mathrm T}    \\
T_{k\tilde{k},nm}=\left(
P^{k11,\tilde{k}00} _{nm}, P^{k10,\tilde{k}01} _{nm}, P^{k01,\tilde{k}10} _{mn}, P^{k00,\tilde{k}11} _{nm}
\right)^{\mathrm T}   \\
\end{array}
\right.
\end{eqnarray}
then
\begin{eqnarray}
\left(
(\langle k,\omega \rangle + \langle \tilde{k},\tilde{\omega} \rangle)I +
A_{n+tc} \otimes I_2 - I_2 \otimes A_{m+tc} ^{\mathrm T}
\right)
Q_{k\tilde{k},n+tc,m+tc} = -i T_{k\tilde{k},n+tc,m+tc}
\end{eqnarray}
By assumption (2.18) in condition $\bf (A3)$:
\begin{eqnarray}
|\det \left(
(\langle k,\omega \rangle + \langle \tilde{k},\tilde{\omega} \rangle)I +
A_{n+tc} \otimes I_2 - I_2 \otimes A_{m+tc} ^{\mathrm T}
\right) | \geq \frac{\gamma}{K^{\tau}}  \nonumber
\end{eqnarray}
Notice that the "main part" of $\Omega_n, \tilde{\Omega}_n$ is $|n|^2$ and
\begin{eqnarray}
\langle k,\omega \rangle + \langle \tilde{k},\tilde{\omega} \rangle + |n+tc|^2 - |m+tc|^2 =
\langle k,\omega \rangle + \langle \tilde{k},\tilde{\omega} \rangle + |n|^2 - |m|^2 +2t\langle n-m,c\rangle
\nonumber
\end{eqnarray}
So if $\langle n-m,c \rangle=0$, then we have
\begin{eqnarray}
\lim_{t\rightarrow \infty} Q_{k\tilde{k},n+tc,m+tc} =
-i\left(
(\langle k,\omega \rangle + \langle \tilde{k},\tilde{\omega} \rangle)I +
\lim_{t\rightarrow \infty}
(
A_{n+tc} \otimes I_2 - I_2 \otimes A_{m+tc} ^{\mathrm T}
)
\right)^{-1}
\lim_{t\rightarrow \infty}T_{k\tilde{k},n+tc,m+tc} \nonumber
\end{eqnarray}
exists. And according to
\begin{eqnarray}
\|
\lim_{t\rightarrow \infty} Q_{k\tilde{k},n+tc,m+tc}
\| \leq
\gamma^{-16} K^{16\tau+16} \varepsilon e^{-(|k|+|\tilde{k}|)r} e^{-|n-m|\rho}
\end{eqnarray}

\begin{eqnarray}
&\left(
(\langle k,\omega \rangle + \langle \tilde{k},\tilde{\omega} \rangle)I +
A_{n+tc} \otimes I_2 - I_2 \otimes A_{m+tc} ^{\mathrm T}
\right)
\left(
Q_{k\tilde{k},n+tc,m+tc} - \lim_{t\rightarrow \infty} Q_{k\tilde{k},n+tc,m+tc}
\right)  \nonumber \\
&  =-i
\left(
T_{k\tilde{k},n+tc,m+tc} - \lim_{t\rightarrow \infty} T_{k\tilde{k},n+tc,m+tc}
\right)   \nonumber \\
&    -
\left(
A_{n+tc} \otimes I_2 - I_2 \otimes A_{m+tc} ^{\mathrm T} -
\lim_{t\rightarrow \infty}
( A_{n+tc} \otimes I_2 - I_2 \otimes A_{m+tc} ^{\mathrm T} )
\right)
\lim_{t\rightarrow \infty} Q_{k\tilde{k},n+tc,m+tc}   \nonumber
\end{eqnarray}
and
\begin{eqnarray}
&\|
\left(
A_{n+tc}\otimes I_2 - I_2 \otimes A_{m+tc} -
\lim_{t\rightarrow \infty} ( A_{n+tc}\otimes I_2 - I_2 \otimes A_{m+tc} )
\right)
\lim_{t\rightarrow \infty}Q_{k\tilde{k},n+tc,m+tc}
\|    \nonumber \\
&\leq
\frac{\varepsilon_0}{|t|}\gamma^{-16} K^{16\tau+16} \varepsilon e^{-(|k|+|\tilde{k}|)r}
e^{-|n-m|\rho}    \nonumber
\end{eqnarray}
So we get
\begin{eqnarray}
\|
Q_{k\tilde{k},n+tc,m+tc} - \lim_{t\rightarrow \infty} Q_{k\tilde{k},n+tc,m+tc}
\|  &\leq&
\gamma^{-32} K^{32\tau+32} \frac{\varepsilon}{|t|} e^{-(|k|+|\tilde{k}|)r} e^{-|n-m|\rho}   \nonumber \\
&\leq& \frac{\varepsilon^{\frac{2}{3}}}{|t|} e^{-(|k|+|\tilde{k}|)r} e^{-|n-m|\rho}
\end{eqnarray}
From (4.61) we conclude that
\begin{eqnarray}
&\|
\frac{\partial^2 F}{\partial u_{n+tc} \partial \bar u_{m+tc} } -
\lim\limits_{t\rightarrow \infty} \frac{\partial^2 F}{\partial u_{n+tc} \partial \bar u_{m+tc} }
\|   \leq
\frac{\varepsilon^{\frac{2}{3}}}{|t|}e^{|n-m|\rho} \nonumber \\
&\|
\frac{\partial^2 F}{\partial u_{n+tc} \partial \bar v_{m+tc} } -
\lim\limits_{t\rightarrow \infty} \frac{\partial^2 F}{\partial u_{n+tc} \partial \bar v_{m+tc} }
\|    \leq
\frac{\varepsilon^{\frac{2}{3}}}{|t|}e^{|n-m|\rho} \nonumber \\
&\|
\frac{\partial^2 F}{\partial v_{n+tc} \partial \bar u_{m+tc} } -
\lim\limits_{t\rightarrow \infty} \frac{\partial^2 F}{\partial v_{n+tc} \partial \bar u_{m+tc} }
\|   \leq
\frac{\varepsilon^{\frac{2}{3}}}{|t|}e^{|n-m|\rho} \nonumber \\
&\|
\frac{\partial^2 F}{\partial v_{n+tc} \partial \bar v_{m+tc} } -
\lim\limits_{t\rightarrow \infty} \frac{\partial^2 F}{\partial v_{n+tc} \partial \bar v_{m+tc} }
\|    \leq
\frac{\varepsilon^{\frac{2}{3}}}{|t|}e^{|n-m|\rho} \nonumber
\end{eqnarray}

If $\langle n-m,c \rangle \neq 0$ and $|t|>K$, we have
\begin{eqnarray}
&\|  \frac{\partial^2 F}{\partial u_{n+tc} \partial \bar u_{m+tc} }  - 0 \| \leq
  \frac{\varepsilon}{|t|} e^{-|n-m|\rho}   \nonumber \\
&\|  \frac{\partial^2 F}{\partial u_{n+tc} \partial \bar v_{m+tc} }  - 0 \| \leq
  \frac{\varepsilon}{|t|} e^{-|n-m|\rho}   \nonumber \\
&\|  \frac{\partial^2 F}{\partial v_{n+tc} \partial \bar u_{m+tc} }  - 0 \| \leq
  \frac{\varepsilon}{|t|} e^{-|n-m|\rho}   \nonumber \\
&\|  \frac{\partial^2 F}{\partial v_{n+tc} \partial \bar v_{m+tc} }  - 0 \| \leq
  \frac{\varepsilon}{|t|} e^{-|n-m|\rho}   \nonumber
\end{eqnarray}

For other related terms, which concern derivatives of $F$ with respect to
$$
u_n u_m, u_n v_m, v_n u_m, v_n v_m;
\bar u_n \bar u_m, \bar u_n \bar v_m, \bar v_n \bar u_m, \bar v_n \bar v_m
$$
similar estimates hold.( In these cases, use
$\frac{\varepsilon}{|t|}e^{-|n+m|\rho}$  instead of $\frac{\varepsilon}{|t|}e^{-|n-m|\rho}$)
And it's easy to see that for the cases when $n\in \Z^d _1 \setminus \Z^d _2, m\in \Z^d _1 \cap \Z^d _2$ and
$n\in \Z^d _1 \cap \Z^d _2, m\in \Z^d _1 \setminus \Z^d _2$ and $n,m \in \Z^d _1 \cap \Z^d _2$ are easier.

To sum up, $F$ satisfies T\"oplitz-Lipschitz property $\bf (A6)$ with $\varepsilon^{\frac{2}{3}}$ in place of $\varepsilon$.  $\Box$

\begin{Lemma}
Assume that $P$ satisfies $\bf (A6)$, $F$ satisfies $\bf (A6)$ with $\varepsilon^{\frac{2}{3}}$ in place of $\varepsilon$ and
\begin{eqnarray}
&\frac{\partial^2 F}{\partial u_n \partial u_m}=0,
 \frac{\partial^2 F}{\partial \bar u_n \partial \bar u_m}=0,
 \frac{\partial^2 F}{\partial v_n \partial v_m}=0,
 \frac{\partial^2 F}{\partial \bar v_n \partial \bar v_m}=0,
\qquad if \quad|n+m| > K  \\
&\frac{\partial^2 F}{\partial u_n \partial \bar u_m}=0,
 \frac{\partial^2 F}{\partial u_n \partial \bar v_m}=0,
 \frac{\partial^2 F}{\partial v_n \partial \bar u_m}=0,
 \frac{\partial^2 F}{\partial v_n \partial \bar v_m}=0,
\qquad if \quad|n-m| > K
\end{eqnarray}
then $\{P,F\}$ satisfies $\bf (A6)$ with $\varepsilon^{\frac{2}{3}}$ in place of $\varepsilon$.
\end{Lemma}

\indent The proof of Lemma 4.5 is the same with Lemma 4.4 in $\cite{12}$, with the help of Lemma 4.4 stated above.

\section{\textbf{Iterative Lemma and Convergence}}
\noindent

For any given $s,\varepsilon,r,\gamma$ and any $\nu \geq 1$, we define the following sequences:
\begin{eqnarray}
&r_{\nu} = r(1-\sum\limits_{i=2} ^{\nu+1} 2^{-i}), \nonumber \\
&\varepsilon_{\nu} = c \gamma^{-16} K_{\nu-1} ^{16\tau+16} \varepsilon_{\nu-1} ^{\frac{4}{3}}  \nonumber \\
&\eta_{\nu}= \varepsilon_{\nu} ^{\frac{1}{3}} \nonumber \\
&L_{\nu}=L_{\nu-1}+\varepsilon_{\nu-1}  \nonumber \\
&s_{\nu}=\frac{1}{4} \eta_{\nu-1} s_{\nu-1}  \nonumber \\
&\rho_{\nu}=\rho (1-\sum\limits_{i=2} ^{\nu+1} 2^{-i})  \nonumber \\
&K_{\nu}=c\left(
(\rho_{\nu-1} - \rho_{\nu})^{-1} \ln \varepsilon_{\nu} ^{-1}
\right)
\end{eqnarray}
where $c$ is a constant, and the parameters $r_0,\varepsilon_0,L_0,s_0,K_0$ are defined to be
$r,\varepsilon,L,s,\ln \frac{1}{\varepsilon}$ respectively.

For later use, we define resonant sets as below:(let ${\mathcal O}_{-1} = {\mathcal O}, K_{-1} = 0 $)
\begin{eqnarray}
{\mathcal R}^{\nu} = \bigcup_{K_{\nu-1} < |k|+|\tilde{k}| \leq K_{\nu}, nm}
 \left(
 {\mathcal R}^{\nu} _{k\tilde{k}} \cup {\mathcal R}^{\nu} _{k\tilde{k},n}
 \cup {\mathcal R}^{\nu,-} _{k\tilde{k},nm} \cup {\mathcal R}^{\nu,+} _{k\tilde{k},nm}
 \right)
\end{eqnarray}
where
\begin{eqnarray}
{\mathcal R}^{\nu} _{k\tilde{k}} &=& \bigg\{
(\xi,\sigma) \in {\mathcal O}_{\nu-1}:
   |\langle k,\omega^{\nu} \rangle + \langle \tilde{k},\tilde{\omega}^{\nu} \rangle |
   < \frac{\gamma}{K_{\nu} ^{\tau}}
\bigg\}  \\
{\mathcal R}^{\nu} _{k\tilde{k},n} &=& \bigg\{
(\xi,\sigma) \in {\mathcal O}_{\nu-1}:
   |\det\left((\langle k,\omega^{\nu} \rangle + \langle \tilde{k},\tilde{\omega}^{\nu} \rangle)I
   +A^{\nu} _n
   \right)|
   < \frac{\gamma}{K_{\nu} ^{\tau}}
\bigg\}  \\
{\mathcal R}^{\nu,-} _{k\tilde{k},nm} &=& \bigg\{
(\xi,\sigma) \in {\mathcal O}_{\nu-1}:
   |\det\bigg((\langle k,\omega^{\nu} \rangle + \langle \tilde{k},\tilde{\omega}^{\nu} \rangle)I  \nonumber \\
   &+&  A^{\nu} _n \otimes I_2 - I_2 \otimes {A^{\nu} _m}^{\mathrm T}
   \bigg)|
   < \frac{\gamma}{K_{\nu} ^{\tau}}
\bigg\}  \\
{\mathcal R}^{\nu,+} _{k\tilde{k},nm} &=& \bigg\{
(\xi,\sigma) \in {\mathcal O}_{\nu-1}:
   |\det\bigg((\langle k,\omega^{\nu} \rangle + \langle \tilde{k},\tilde{\omega}^{\nu} \rangle)I \nonumber \\
   &+&   A^{\nu} _n \otimes I_2 + I_2 \otimes {A^{\nu} _m}
   \bigg)|
   < \frac{\gamma}{K_{\nu} ^{\tau}}
\bigg\}
\end{eqnarray}

\begin{Lemma}(Iterative Lemma)

Let $\varepsilon > 0$ be sufficiently small, $\nu\geq 0$. Suppose that

\noindent (1)
$
N_{\nu}+{\mathcal B}_{\nu} = \langle \omega,I \rangle+ \langle \tilde{\omega},J \rangle +
\sum\limits_{n\in \Z^d _1} \Omega_n u_n \bar u_n + \sum\limits_{n\in \Z^d _2} \tilde{\Omega}_n v_n \bar v_n +
\sum\limits_{n\in \Z^d _1 \cap \Z^d _2} (a_n u_n \bar v_n + b_n \bar u_n v_n)
$
is a normal form with parameters $(\xi,\sigma)$ in a closed set ${\mathcal O}_{\nu}$ of $\R^{b+\tilde{b}}$ satisfying
\begin{eqnarray}
&|\langle k,\omega_{\nu} \rangle+\langle \tilde{k},\tilde{\omega}_{\nu} \rangle| \geq
      \frac{\gamma}{K_{\nu} ^{\tau}},   \quad
      0< |k|+|\tilde{k}|\leq K_{\nu}   \nonumber \\
&|\det\left(
(\langle k,\omega_{\nu} \rangle+\langle \tilde{k},\tilde{\omega}_{\nu} \rangle)I + A_n ^{\nu}
\right)|  \geq \frac{\gamma}{K_{\nu} ^{\tau}}, \quad
     |k|+|\tilde{k}| \leq K_{\nu}  \nonumber  \\
&|\det\left(
(\langle k,\omega_{\nu} \rangle+\langle \tilde{k},\tilde{\omega}_{\nu} \rangle) I +
A^{\nu} _n \otimes I_2 - I_2 \otimes {A^{\nu} _m}^{\mathrm T}
\right)|  \geq \frac{\gamma}{K_{\nu} ^{\tau}}, \quad
      |k|+|\tilde{k}|\leq K_{\nu}, |k|+|\tilde{k}| +||n|-|m|| \neq 0  \nonumber \\
&|\det\left(
(\langle k,\omega_{\nu} \rangle+\langle \tilde{k},\tilde{\omega}_{\nu} \rangle) I +
A^{\nu} _n \otimes I_2 + I_2 \otimes {A^{\nu}} _m
\right)|  \geq \frac{\gamma}{K_{\nu} ^{\tau}}, \quad
      |k|+|\tilde{k}|\leq K_{\nu}  \nonumber
\end{eqnarray}
where
\begin{eqnarray}
A_n =
\left(
\begin{array}{cc}
\Omega_n  & a_n \\
b_n & \ \tilde{\Omega}_n
\end{array}
\right)  \qquad n \in \Z^d _1 \cap \Z^d _2  \nonumber
\end{eqnarray}
and
\begin{eqnarray}
A_n = \Omega_n, \quad n \in \Z^d _1 \setminus \Z^d _2   \nonumber \\
A_n = \tilde{\Omega}_n, \quad n \in \Z^d _2 \setminus \Z^d _1   \nonumber
\end{eqnarray}

\noindent (2)
$
\omega_{\nu}(\xi,\sigma),\tilde{\omega}_{\nu}(\xi,\sigma),
\Omega_{\nu}(\xi,\sigma),\tilde{\Omega}_{\nu}(\xi,\sigma)
$
are $C^4 _W$ smooth in $(\xi,\sigma)$ satisfying
\begin{eqnarray}
&|\omega_{\nu} - \omega_{\nu-1}|_{{\mathcal O}_{\nu}}  < \varepsilon_{\nu-1},
|\tilde{\omega}_{\nu} - \tilde{\omega}_{\nu-1}|_{{\mathcal O}_{\nu}}  < \varepsilon_{\nu-1},
|\Omega_{\nu} - \Omega_{\nu-1}|_{{\mathcal O}_{\nu}}  < \varepsilon_{\nu-1},
|\tilde{\Omega}_{\nu} - \tilde{\Omega}_{\nu-1}|_{{\mathcal O}_{\nu}}  < \varepsilon_{\nu-1}  \nonumber \\
&|a_n ^{\nu}-a_n ^{\nu-1}|_{{\mathcal O}_{\nu}} < \varepsilon_{\nu-1},
 |b_n ^{\nu}-b_n ^{\nu-1}|_{{\mathcal O}_{\nu}} < \varepsilon_{\nu-1}  \nonumber
\end{eqnarray}

\noindent (3)
$N_{\nu} + {\mathcal B}_{\nu} + P_{\nu}$ satisfies $\bf (A5),(A6)$ with parameters $K_{\nu},\varepsilon_{\nu},\rho_{\nu}$ and
$$
\| X_{P_{\nu}} \|_{D(r_{\nu}, s_{\nu}), {\mathcal O}_{\nu}} < \varepsilon_{\nu}
$$
Then there is a closed subset ${\mathcal O}_{\nu+1} \subseteq {\mathcal O}_{\nu}$ with
$$
{\mathcal O}_{\nu+1} = {\mathcal O}_{\nu} \setminus {\mathcal R}^{\nu+1}
$$
where ${\mathcal R}^{\nu+1}$ is defined in (5.2). We have a symplectic transformation of variables:
$$
\Phi_{\nu} : D_{\rho_{\nu+1}} (r_{\nu+1},s_{\nu+1}) \times {\mathcal O}_{\nu} \rightarrow
D_{\rho_{\nu}} (r_{\nu},s_{\nu})
$$
s.t. in
$
D_{\rho_{\nu+1}} (r_{\nu+1},s_{\nu+1}) \times {\mathcal O}_{\nu}
$,
$H_{\nu+1} = H_{\nu} \circ \Phi_{\nu}$ has the form:
\begin{eqnarray}
H_{\nu+1}=\langle k,\omega_{\nu+1} \rangle + \langle \tilde{k},\tilde{\omega}_{\nu+1} \rangle +
          \sum_{n\in \Z^d _1} \Omega^{\nu+1} _n u_n \bar u_n +
          \sum_{n\in \Z^d _2} \Omega^{\nu+1} _n v_n \bar v_n +
          \sum_{n\in \Z^d _1 \cap \Z^d _2} (a^{\nu+1} _n u_n \bar v_n + b^{\nu+1} _n \bar u_n v_n )  \nonumber
\end{eqnarray}
with
\begin{eqnarray}
&|\omega_{\nu+1} - \omega_{\nu}|_{{\mathcal O}_{\nu+1}}  < \varepsilon_{\nu},
|\tilde{\omega}_{\nu+1} - \tilde{\omega}_{\nu}|_{{\mathcal O}_{\nu+1}}  < \varepsilon_{\nu},
|\Omega_{\nu+1} - \Omega_{\nu}|_{{\mathcal O}_{\nu+1}}  < \varepsilon_{\nu},
|\tilde{\Omega}_{\nu+1} - \tilde{\Omega}_{\nu}|_{{\mathcal O}_{\nu+1}}  < \varepsilon_{\nu}  \nonumber \\
&|a_n ^{\nu+1}-a_n ^{\nu}|_{{\mathcal O}_{\nu+1}} < \varepsilon_{\nu},
 |b_n ^{\nu+1}-b_n ^{\nu}|_{{\mathcal O}_{\nu+1}} < \varepsilon_{\nu}  \nonumber
\end{eqnarray}

For parameter $(\xi,\eta)$ in ${\mathcal O}_{\nu+1}$, we have the following Diophantine condition:
\begin{eqnarray}
&|\langle k,\omega_{\nu+1} \rangle+\langle \tilde{k},\tilde{\omega}_{\nu+1} \rangle| \geq
      \frac{\gamma}{K_{\nu+1} ^{\tau}},   \quad
      0< |k|+|\tilde{k}|\leq K_{\nu+1}   \nonumber \\
&|\det\left(
(\langle k,\omega_{\nu+1} \rangle+\langle \tilde{k},\tilde{\omega}_{\nu+1} \rangle)I + A_n ^{\nu+1}
\right)|  \geq \frac{\gamma}{K_{\nu+1} ^{\tau}}, \quad
     |k|+|\tilde{k}| \leq K_{\nu+1}  \nonumber  \\
&|\det\left(
(\langle k,\omega_{\nu+1} \rangle+\langle \tilde{k},\tilde{\omega}_{\nu+1} \rangle) I +
A^{\nu+1} _n \otimes I_2 - I_2 \otimes {A^{\nu+1} _m}^{\mathrm T}
\right)| \geq   \frac{\gamma}{K_{\nu+1} ^{\tau}}, \nonumber \\
    &     |k|+|\tilde{k}|\leq K_{\nu+1}, |k|+|\tilde{k}| +||n|-|m|| \neq 0  \nonumber \\
&|\det\left(
(\langle k,\omega_{\nu+1} \rangle+\langle \tilde{k},\tilde{\omega}_{\nu+1} \rangle) I +
A^{\nu+1} _n \otimes I_2 + I_2 \otimes {A^{\nu+1}} _m \right)| \geq  \frac{\gamma}{K_{\nu+1} ^{\tau}},
    \nonumber \\
  &   |k|+|\tilde{k}|\leq K_{\nu+1}  \nonumber
\end{eqnarray}

Besides, $N_{\nu+1} + {\mathcal B}_{\nu+1} + P_{\nu+1}$ also satisfies condition $\bf (A5),(A6)$ with
$K_{\nu+1}, \varepsilon_{\nu+1}, \rho_{\nu+1}$ in place of $K_{\nu}, \varepsilon_{\nu}, \rho_{\nu}$
and
$$
\|X_{P_{\nu+1}}\|_{D_{\rho_{\nu+1}} (r_{\nu+1}, s_{\nu+1}), {\mathcal O}_{\nu+1}} \leq \varepsilon_{\nu+1}
$$
\end{Lemma}

(Convergence):
Suppose that the assumption of Theorem 2 are satisfied. Recall that
\begin{eqnarray}
\varepsilon_0 = \varepsilon,
r_0 = r,
s_0 = s,
\rho_0 = \rho,
L_0 = L,
N_0 = N,
{\mathcal B}_0 = 0,
P_0 = P  \nonumber
\end{eqnarray}
and $\mathcal O$ is a bounded positive-measure set in $\R^{b+\tilde{b}}$. The assumptions of the iteration lemma are satisfied when $\nu=0$ if $\varepsilon_0$ and $\gamma$ are sufficiently small. Inductively, we obtain the following sequences:
\begin{eqnarray}
& {\mathcal O}_{\nu+1} \subseteq {\mathcal O}_{\nu}, \nonumber \\
& \Psi^{\nu} = \Phi_0 \circ \Phi_1 \circ \cdots \circ \Phi_{\nu}:
   D_{\rho_{\nu+1}} (r_{\nu+1},s_{\nu+1}) \times {\mathcal O}_{\nu} \rightarrow D_{\rho_0}(r_0,s_0), \nu \geq 0
   \nonumber \\
& H \circ \Psi^{\nu}= H_{\nu+1}=N_{\nu+1} + {\mathcal B}_{\nu+1} +  P_{\nu+1}  \nonumber
\end{eqnarray}

Let $\tilde{{\mathcal O}}=\bigcap_{\nu=0} ^{\infty} {\mathcal O}_{\nu}$. By the help of Lemma 4.2, we conclude that
$
N_{\nu}, \Psi^{\nu}, D\Psi^{\nu}, \omega_{\nu}, \tilde{\omega}_{\nu}
$
converge uniformly on $D_{\frac{1}{2}\rho} (\frac{1}{2}r,0) \times \tilde{{\mathcal O}}$ with
\begin{eqnarray}
N_{\infty} + {\mathcal B}_{\infty} =
\langle \omega_{\infty},I \rangle + \langle \tilde{\omega}_{\infty},J \rangle +
\sum_{n\in \Z^d _1} \Omega^{\infty} _n u_n \bar u_n +
\sum_{n\in \Z^d _2} \tilde{\Omega}^{\infty} _n v_n \bar v_n +
\sum_{n\in \Z^d _1 \cap \Z^d _2} (a^{\infty} _n u_n \bar v_n + b^{\infty} _n \bar u_n v_n)  \nonumber
\end{eqnarray}
Since $\varepsilon_{\nu+1} = c\gamma^{-16}K^{16\tau+16} \varepsilon_{\nu} ^{\frac{4}{3}}$, it follows that
$\varepsilon_{\nu+1} \rightarrow 0$ provided that $\varepsilon$ is small enough.

Let $\phi^t _H$ be the flow of $X_H$. Since $H\circ \Psi^{\nu} = H_{\nu+1}$, we have
\begin{eqnarray}
\phi^t _H \circ \Psi^{\nu} = \Psi^{\nu} \circ \phi^t _{H_{\nu+1}}
\end{eqnarray}
The uniform convergence of $\Psi^{\nu}, D\Psi^{\nu}, \omega_{\nu}, \tilde{\omega}_{\nu}, X_{H_{\nu}}$ implies that the limits can be taken on both sides of (5.7). Hence on
$
D_{\frac{1}{2}\rho}(\frac{1}{2}r,0) \times \tilde{{\mathcal O}}
$
we get
\begin{eqnarray}
\phi^t _H \circ \Psi^{\infty} = \Psi^{\infty} \circ \phi^t _{H_{\infty}}
\end{eqnarray}
and
$$
\Psi^{\infty}:D_{\frac{1}{2}\rho}(\frac{1}{2}r,0) \times \tilde{{\mathcal O}} \rightarrow
D_{\rho}(r,s) \times {\mathcal O}
$$
It follows from (5.8) that
$$
\phi^t _H \left(
\Psi^{\infty} (\T^{b+\tilde{b}} \times \{\xi,\sigma\})
\right)  =
\Psi^{\infty} (\T^{b+\tilde{b}} \times \{\xi,\sigma\})
$$
for $(\xi,\sigma) \in \tilde{{\mathcal O}}$. This means that
$\Psi^{\infty} (\T^{b+\tilde{b}} \times \{\xi,\sigma\} )$ is an embedded torus which is invariant for the original perturbed Hamiltonian system at $(\xi,\sigma) \in \tilde{{\mathcal O}}$. We remark that the frequencies
$(\omega^{\infty}(\xi,\sigma),\tilde{\omega}^{\infty}(\xi,\sigma))$ associated to
$\Psi^{\infty} (\T^{b+\tilde{b}} \times \{\xi+\sigma\})$ are slightly different from the unperturbed frequencies
$(\omega(\xi,\sigma),\tilde{\omega}(\xi,\sigma))$.

\section{\textbf{Measure Estimate}}
\noindent

Recall the definition of resonant sets (5.2)-(5.6). For convenience, we only deal with the most complicated case (5.5), and we assume $n,m \in \Z^d _1 \cap \Z^d _2$(when $n\in \Z^d _1 \setminus \Z^d _2$ or
$n\in \Z^d _1 \setminus \Z^d _2$, the proof is easier). We have known that at the $\nu$-th step, the small divisor conditions are automatically satisfied for $|k|+|\tilde{k}| \leq K_{\nu}$, so we only need to verify the case when $K_{\nu}< |k|+|\tilde{k}|\leq K_{\nu+1}$. And for simplicity, we just set
\begin{eqnarray}
M^{\nu}=(\langle \omega_{\nu},I \rangle + \langle \tilde{\omega}_{\nu},J \rangle)I + A_n ^{\nu} \otimes I_2 -
     I_2 \otimes {A_m ^{\nu}}^{\mathrm T}
\end{eqnarray}

To make use of condition $\bf (A6)$: T\"oplitz-Lipschitz property, we state two technical lemmas here:
\begin{Lemma}(Lemma 6.1 in $\cite{17}$)

For any given $n,m \in \Z^d _1 \cap \Z^d _2$ with $|n-m|\leq K_{\nu}$, either
$
|\det\big(
(\langle \omega_{\nu},I \rangle + \langle \tilde{\omega}_{\nu},J \rangle)I +
A_n ^{\nu} \otimes I_2 - I_2 \otimes {A_m ^{\nu}}^{\mathrm T}
\big)| \geq 1
$
or there are
$
n_0,m_0,c_1,\cdots,c_{d-1} \in \Z^d
$
with
$
|n_0|,|m_0|,|c_1|,\cdots,|c_{d-1}| \leq 3K_{\nu} ^2
$
and
$t_1,t_2,\cdots,t_{d-1} \in Z$, such that
$
n=n_0 + \sum\limits_{j=1} ^{d-1} t_j c_j,
m=m_0 + \sum\limits_{j=1} ^{d-1} t_j c_j
$.
\end{Lemma}

We omit its proof.

\begin{Lemma}(Lemma 6.2 in $\cite{17}$)

$$
\bigcup_{k\tilde{k};n,m \in \Z^d _1 \cap \Z^d _2} {\mathcal R}^{\nu,-} _{k\tilde{k},nm} \subseteq
\bigcup_{n_0,m_0,c_1,\cdots,c_{d-1} \in \Z^d;t_1,\cdots,t_{d-1}\in \Z}
{\mathcal R}^{\nu,-} _
{
k\tilde{k},n_0 + \sum_{j=1} ^{d-1} t_j c_j , m_0 + \sum_{j=1} ^{d-1} t_j c_j
}
$$
where
$
|n_0|,|m_0|,|c_1|,\cdots,|c_{d-1}| \leq 3K_{\nu} ^2
$.
\end{Lemma}

We omit its proof.

\begin{Lemma}
For fixed $k,\tilde{k},n_0,m_0,c_1,c_2,\cdots,c_{d-1}$, one has
$$
meas\left(
\bigcup_{t_1,\cdots,t_{d-1} \in \Z} {\mathcal R}^{\nu,-} _
{
k\tilde{k},n_0 + \sum_{j=1} ^{d-1} t_j c_j,m_0 + \sum_{j=1} ^{d-1} t_j c_j
}
\right)  < c\frac{\gamma}{K_{\nu} ^{\frac{\tau}{d!}}}
$$
\end{Lemma}

\noindent $\bf Proof:$
Without loss of generality, we could just assume that
$
|t_1|\leq |t_2| \leq \cdots \leq |t_{d-1}|
$.
Recall that
$
\Omega_n = |n|^2 + \acute{\Omega}_n, \tilde{\Omega}_n = |n|^2 + \acute{\tilde{\Omega}}_n
$
and according to condition $\bf (A6)$, we have
\begin{eqnarray}
&|
{\acute{\Omega}}^{\nu} _{n_0 + \sum_{j=1} ^{d-1} t_j c_j }  -
   \lim\limits_{t_j \rightarrow \infty} {\acute{\Omega}}^{\nu} _{n_0 + \sum_{j=1} ^{d-1} t_j c_j }
|   <  \frac{\varepsilon_0}{|t_j|}, 1\leq j \leq d-1   \nonumber \\
&|
{\acute{\Omega}}^{\nu} _{m_0 + \sum_{j=1} ^{d-1} t_j c_j }  -
   \lim\limits_{t_j \rightarrow \infty} {\acute{\Omega}}^{\nu} _{m_0 + \sum_{j=1} ^{d-1} t_j c_j }
|   <  \frac{\varepsilon_0}{|t_j|}, 1\leq j \leq d-1   \nonumber \\
&|
{\acute{\tilde{\Omega}}}^{\nu} _{n_0 + \sum_{j=1} ^{d-1} t_j c_j }  -
   \lim\limits_{t_j \rightarrow \infty} {\acute{\tilde{\Omega}}}^{\nu} _{n_0 + \sum_{j=1} ^{d-1} t_j c_j }
|   <  \frac{\varepsilon_0}{|t_j|}, 1\leq j \leq d-1   \nonumber \\
&|
{\acute{\tilde{\Omega}}}^{\nu} _{m_0 + \sum_{j=1} ^{d-1} t_j c_j }  -
   \lim\limits_{t_j \rightarrow \infty} {\acute{\tilde{\Omega}}}^{\nu} _{m_0 + \sum_{j=1} ^{d-1} t_j c_j }
|   <  \frac{\varepsilon_0}{|t_j|}, 1\leq j \leq d-1   \nonumber \\
&|
 a^{\nu} _{n_0 + \sum_{j=1} ^{d-1} t_j c_j }  -
   \lim\limits_{t_j \rightarrow \infty} a^{\nu} _{n_0 + \sum_{j=1} ^{d-1} t_j c_j }
|   <  \frac{\varepsilon_0}{|t_j|}, 1\leq j \leq d-1   \nonumber \\
&|
 a^{\nu} _{m_0 + \sum_{j=1} ^{d-1} t_j c_j }  -
   \lim\limits_{t_j \rightarrow \infty} a^{\nu} _{m_0 + \sum_{j=1} ^{d-1} t_j c_j }
|   <  \frac{\varepsilon_0}{|t_j|}, 1\leq j \leq d-1   \nonumber \\
&|
 b^{\nu} _{n_0 + \sum_{j=1} ^{d-1} t_j c_j }  -
   \lim\limits_{t_j \rightarrow \infty} b^{\nu} _{n_0 + \sum_{j=1} ^{d-1} t_j c_j }
|   <  \frac{\varepsilon_0}{|t_j|}, 1\leq j \leq d-1   \nonumber \\
&|
 b^{\nu} _{m_0 + \sum_{j=1} ^{d-1} t_j c_j }  -
   \lim\limits_{t_j \rightarrow \infty} b^{\nu} _{m_0 + \sum_{j=1} ^{d-1} t_j c_j }
|   <  \frac{\varepsilon_0}{|t_j|}, 1\leq j \leq d-1   \nonumber
\end{eqnarray}

Hence we have
\begin{eqnarray}
|
\det\bigg(M^{\nu}(t_1,\cdots,t_{d-1})\bigg) -
\lim\limits_{t_j \rightarrow \infty} \det\bigg(M^{\nu}(t_1,\cdots,t_{d-1})\bigg)
|  < \frac{\varepsilon_0}{|t_j|}, 1\leq j \leq d-1  \nonumber
\end{eqnarray}

Now we introduce the resonant set:
\begin{eqnarray}
{\mathcal R}^{\nu} _{k\tilde{k},n_0 m_0 c_1 \cdots c_{d-1} \infty^{d-1}} =
      \{
   (\xi,\sigma)\in {\mathcal O}_{\nu-1}:
|
\lim\limits_{t_1 \rightarrow \infty}\bigg(
\lim_{t_2,\cdots,t_{d-1}\rightarrow \infty}M^{\nu} (t_1,\cdots,t_{d-1})
          \bigg)
|  <  \frac{\gamma}{K_{\nu} ^{\frac{\tau}{d!}}}
      \}
\end{eqnarray}

For fixed $k,\tilde{k},n_0,m_0,c_1,\cdots,c_{d-1}$, we have the estimate
\begin{eqnarray}
meas\left(
{\mathcal R}^{\nu} _{k\tilde{k},n_0m_0c_1 \cdots c_{d-1} \infty^{d-1}}
\right) < \frac{\gamma}{K_{\nu} ^{\frac{\tau}{d!}}}
\end{eqnarray}

So for
$
(\xi,\sigma) \in {\mathcal O}_{\nu-1} \setminus
{\mathcal R}^{\nu} _{k\tilde{k},n_0m_0c_1 \cdots c_{d-1} \infty^{d-1}}
$, we have
\begin{eqnarray}
|
\lim_{t_1\rightarrow\infty}
\left(
\lim_{t_2,\cdots,t_{d-1} \rightarrow \infty}M^{\nu} (t_1,\cdots,t_{d-1})
\right)
| \geq \frac{\gamma}{K_{\nu} ^{\frac{\tau}{d!}}}
\end{eqnarray}

\indent Now we consider the following $d$ kinds of cases, according the size of each $t_j$:

\indent Case $\bf (1)$:
If $|t_1|>K_{\nu} ^{\frac{\tau}{d!}}$, then for
$
(\xi,\sigma) \in {\mathcal O}_{\nu-1} \setminus
{\mathcal R}^{\nu} _{k\tilde{k},n_0m_0c_1 \cdots c_{d-1} \infty^{d-1}}
$, we have
\begin{eqnarray}
&&|M^{\nu} (t_1,\cdots,t_{d-1})| \nonumber \\
&&\geq |\lim_{t_2,\cdots,t_{d-1}}\rightarrow\infty M^{\nu} (t_1,\cdots,t_{d-1}) |
   -  \sum_{j=1} ^{d-1} \frac{\varepsilon_0}{|t_j|}  \nonumber \\
&&\geq \frac{\gamma}{K_{\nu} ^{\frac{\tau}{d!}}}
    -(d-1) \frac{\varepsilon_0}{K_{\nu} ^{\frac{\tau}{d!}}} \nonumber \\
&&\geq \frac{\gamma}{2K_{\nu} ^{\frac{\tau}{d!}}}
\end{eqnarray}

\indent Case $\bf (2)$: When $|t_1|\leq K_{\nu} ^{\frac{\tau}{d!}},|t_2|\geq K_{\nu} ^{\frac{2\tau}{d!}}$, we define the resonant:
\begin{eqnarray}
{\mathcal R}^{\nu} _{k\tilde{k},n_0m_0c_1 \cdots c_{d-1} t_1 \infty^{d-2}} =
\{
(\xi,\sigma)\in {\mathcal O}_{\nu-1}:
|
\lim_{t_2,\cdots,t_{d-1}\rightarrow\infty}
  M^{\nu} (t_1,\cdots,t_{d-1})
| < \frac{\gamma}{K_{\nu} ^{\frac{2\tau}{d!}}}
\}
\end{eqnarray}

For fixed $k,\tilde{k},n_0,m_0,c_1,\cdots,c_{d-1},t_1$, one has
\begin{eqnarray}
meas(
{\mathcal R}^{\nu} _
   {
k\tilde{k},n_0m_0c_1\cdots c_{d-1} t_1 \infty^{d-2}
   }
) < \frac{\gamma}{K_{\nu} ^{\frac{2\tau}{d!}}}
\end{eqnarray}

then

\begin{eqnarray}
meas\{
\cup_{|t_1|\leq K_{\nu} ^{\frac{\tau}{d!}}}
{\mathcal R}^{\nu} _{k\tilde{k},n_0m_0c_1\cdots c_{d-1}t_1 \infty^{d-2}}
\} < K_{\nu} ^{\frac{\tau}{d!}} \times \frac{\gamma}{K_{\nu} ^{\frac{2\tau}{d!}}}
  < \frac{\gamma}{K_{\nu} ^{\frac{\tau}{d!}}}
\end{eqnarray}

So if
$
|t_1|\leq K_{\nu} ^{\frac{\tau}{d!}}, |t_2|\geq K_{\nu} ^{\frac{2\tau}{d!}}
$,
for
$
(\xi,\sigma)\in {\mathcal O}_{\nu-1} \setminus {\mathcal R}^{\nu} _
{
k\tilde{k},n_0m_0c_1 \cdots c_{d-1} t_1 \infty^{d-2}
}
$, we have
\begin{eqnarray}
&&|M^{\nu} (t_1,\cdots,t_{d-1})| \nonumber \\
&&\geq |\lim_{t_2,\cdots,t_{d-1}}\rightarrow\infty M^{\nu} (t_1,\cdots,t_{d-1}) |
   -  \sum_{j=2} ^{d-1} \frac{\varepsilon_0}{|t_j|}  \nonumber \\
&&\geq \frac{\gamma}{K_{\nu} ^{\frac{2\tau}{d!}}} -
   (d-2) \frac{\varepsilon_0}{K_{\nu} ^{\frac{2\tau}{d!}}} \nonumber \\
&&\geq \frac{\gamma}{2K_{\nu} ^{\frac{2\tau}{d!}}}
\end{eqnarray}

We can continue this process until the $d-1$ th step:

\indent Case $\bf (d-1)$:
If
$
|t_1|\leq K_{\nu} ^{\frac{\tau}{d!}},|t_2|\leq K_{\nu} ^{\frac{2\tau}{d!}},\cdots,
|t_{d-2}| \leq K_{\nu} ^{\frac{(d-2)!\tau}{d!}},
|t_{d-1}| \geq K_{\nu} ^{\frac{(d-1)!\tau}{d!}}
$, we define the resonant set

\begin{eqnarray}
{\mathcal R}^{\nu} _{k\tilde{k},n_0m_0c_1\cdots c_{d-1},t_1 \cdots t_{d-2} \infty} =
\{
(\xi,\sigma) \in {\mathcal O}_{\nu-1} :
|
\lim\limits_{t_{d-1}\rightarrow\infty} M^{\nu}(t_1,\cdots,t_{d-1})
| < \frac{\gamma}{K_{\nu} ^{\frac{(d-1)!\tau}{d!}}}
\}
\end{eqnarray}

For fixed $k,\tilde{k},n_0,m_0,c_1,\cdots,c_{d-1},t_1,\cdots,t_{d-2}$, one has
\begin{eqnarray}
meas(
{\mathcal R}^{\nu} _
{
k\tilde{k},n_0m_0c_1\cdots c_{d-1},t_1\cdots t_{d-2} \infty
}
) < \frac{\gamma}{K_{\nu} ^{\frac{(d-1)!\tau}{d!}}}
\end{eqnarray}

so we have
\begin{eqnarray}
&meas \left\{
\bigcup_{|t_1|,\cdots,|t_{d-2}|\leq K_{\nu} ^{\frac{(d-2)!\tau}{d!}}}
   {\mathcal R}^{\nu} _
{
k\tilde{k},n_0m_0c_1 \cdots c_{d-1} t_1 \cdots t_{d-2} \infty
}
\right\}    \nonumber \\
&  < K_{\nu} ^{\frac{(d-2)!\times (d-2)\tau}{d!}} \times
   \frac{\gamma}{K_{\nu} ^{\frac{(d-1)!\tau}{d!}}}  \leq \frac{\gamma}{K_{\nu} ^{\frac{\tau}{d!}}}
\end{eqnarray}

Hence if
$
|t_1| \leq K_{\nu} ^{\frac{\tau}{d!}},
|t_2| \leq K_{\nu} ^{\frac{2\tau}{d!}},  \cdots,
|t_{d-2}| \leq K_{\nu} ^{\frac{(d-2)!\tau}{d!}},
|t_{d-1}| > K_{\nu} ^{\frac{(d-1)!\tau}{d!}},
$,
for
$
(\xi,\sigma) \in {\mathcal O}_{\nu-1} \setminus
{\mathcal R}^{\nu} _
{
k\tilde{k},n_0 m_0 c_1 \cdots c_{d-1} t_1 \cdots t_{d-2} \infty
}
$,
we have
\begin{eqnarray}
&&|M^{\nu} (t_1,\cdots, t_{d-1})|  \nonumber \\
&&\geq |\lim\limits_{t_{d-1} \rightarrow \infty} M^{\nu} (t_1,\cdots,t_{d-1})| - \frac{\varepsilon_0}{|t_{d-1}|}
      \nonumber \\
&&\geq \frac{\gamma}{K_{\nu} ^{\frac{(d-1)!\tau}{d!}}} - \frac{\varepsilon_0}{K_{\nu} ^{\frac{(d-1)!\tau}{d!}}}
      \nonumber \\
&&\geq \frac{\gamma}{2K_{\nu} ^{\frac{(d-1)!\tau}{d!}}}
\end{eqnarray}

\indent Case ($\bf d$):
If
$
|t_1| \leq K_{\nu} ^{\frac{\tau}{d!}}, |t_2| \leq K_{\nu} ^{\frac{2\tau}{d!}},\cdots,
|t_{d-1}| \leq K_{\nu} ^{\frac{(d-1)!\tau}{d!}},
$, we define the resonant set
\begin{eqnarray}
{\mathcal R}^{\nu} _{k\tilde{k},n_0m_0c_1\cdots c_{d-1} t_1\cdots t_{d-1}} =
\{
(\xi,\sigma) \in {\mathcal O}_{\nu-1} :
|M^{\nu} (t_1,\cdots,t_{d-1})| < \frac{\gamma}{K_{\nu} ^{\tau}}
\}
\end{eqnarray}
For fixed
$
k,\tilde{k};n_0,m_0,c_1,\cdots, c_{d-1}; t_1,\cdots, t_{d-1}
$,
we have
\begin{eqnarray}
meas\left(
{\mathcal R}^{\nu} _{k\tilde{k},n_0m_0c_1\cdots c_{d-1},t_1\cdots t_{d-1}}
\right) < \frac{\gamma}{K_{\nu} ^{\tau}}
\end{eqnarray}
So we have
\begin{eqnarray}
&& meas\bigg\{
\bigcup_{|t_1|,\cdots,|t_{d-1}|\leq K_{\nu} ^{\frac{(d-1)!\tau}{d!}}}
{\mathcal R}^{\nu} _
{
k\tilde{k},n_0m_0c_1\cdots c_{d-1} t_1\cdots t_{d-1}
}
\bigg\}  \nonumber \\
&& < K_{\nu} ^{\frac{(d-1)!\times (d-1)\tau}{d!}} \times \frac{\gamma}{K_{\nu} ^{\tau}}
    < \frac{\gamma}{K_{\nu} ^{\frac{\tau}{d}}}
\end{eqnarray}

As a consequence, we conclude that
\begin{eqnarray}
meas\left(
\bigcup_{t_1,\cdots,t_{d-1} \in \Z}  {\mathcal R}^{\nu,-} _
{
k\tilde{k},n_0+ \sum\limits_{j=1} ^{d-1} t_j c_j , m_0+ \sum\limits_{j=1} ^{d-1} t_j c_j
}
\right)   <  c\frac{\gamma}{K_{\nu} ^{\frac{\tau}{d!}}}
\end{eqnarray}

So we proved Lemma 6.1. \quad $\Box$

$\\$
\indent By Lemma 6.1, we could get estimates as below:
\begin{eqnarray}
& meas\bigg(
\bigcup_{K_{\nu-1}\leq |k|+|\tilde{k}| \leq K_{\nu}}
     {\mathcal R}^{\nu} _{k\tilde{k}}
  \bigg)  < c\frac{\gamma}{K_{\nu} ^{\tau-b}}
      \nonumber \\
&  meas\bigg(
\bigcup_{K_{\nu-1}\leq |k|+|\tilde{k}| \leq K_{\nu},n}
     {\mathcal R}^{\nu} _{k\tilde{k},n}
  \bigg)  < c\frac{\gamma}{K_{\nu} ^{\tau-d-b}}
      \nonumber \\
&  meas\bigg(
\bigcup_{K_{\nu-1}\leq |k|+|\tilde{k}| \leq K_{\nu},n,m}
     {\mathcal R}^{\nu,\pm} _{k\tilde{k},n,m}
  \bigg)  < c\frac{\gamma}{K_{\nu} ^{\frac{\tau}{d!}-2d(d+1)-(b+\tilde{b})}}
      \nonumber
\end{eqnarray}

Just let $\tau>d!(2d(d+1)+b+\tilde{b}+1)$, then the total measure we need to exclude along the KAM iteration is
\begin{eqnarray}
  meas(\bigcup_{\nu\geq 0} {\mathcal R}^{\nu}) =
  meas\left(
\bigcup_{\nu \geq 0}\bigg(
  \bigcup_{K_{\nu-1}\leq |k|+|\tilde{k}| \leq K_{\nu},n,m}
   {\mathcal R}^{\nu} _k \cup {\mathcal R}^{\nu} _{k,n} \cup {\mathcal R}^{\nu,\pm} _{k,n,m}
                    \bigg)
  \right)  < c\gamma  \nonumber
\end{eqnarray}

So we finished the measure estimate.

\thispagestyle{empty}

\end{document}